\documentclass[12pt]{amsart}
\usepackage{amsmath,amssymb,amscd,array}

\usepackage{hyperref,compress}

\usepackage{DLdef1}

\newcommand {\fby}{{\mathfrak{by}}}

\newcommand {\fj}{{\mathfrak{j}}}
\newcommand {\fme}{{\mathfrak{me}}}

\newcommand{\uN}{\underline{N}}
\newcommand{\del}{\partial}

\begin{document}

\title{Lie algebra deformations in characteristic 2}

\author{Sofiane Bouarroudj${}^1$, Alexei Lebedev${}^2$, Dimitry Leites${}^3$, Irina
Shchepochkina${}^4$}

\address{${}^1$New York University Abu Dhabi\\
Division of Science and Mathematics\\ P.O. Box 129188, United Arab
Emirates; sofiane.bouarroudj@nyu.edu\\
${}^2$Equa
Simulation AB\\
Stockholm, Sweden; alexeylalexeyl@mail.ru\\
${}^3$Department of Mathematics\\
Stockholm University\\
Roslagsv. 101, Kr\"aft\-riket hus 6, SE-106 91 Stockholm, Sweden;
mleites@math.su.se\\
${}^4$Independent University of Moscow\\
Bolshoj Vlasievsky per, dom 11, RU-119 002 Mos\-cow,
Russia;irina@mccme.ru}


%\@namedef{subjclassname@2010}{%
%  \textup{2010} Mathematics Subject Classification}


\begin{abstract}
Of four types of Kaplansky algebras, type-2 and type-4 algebras have
previously unobserved $\mathbb{Z}/2$-gra\-dings: nonlinear in
roots. A method assigning a simple Lie superalgebra to every
$\mathbb{Z}/2$-graded simple Lie algebra in characteristic 2 is
illustrated by seven new series. Type-2 algebras and one of the two type-4
algebras are demystified as nontrivial deforms (the
results of deformations) of the alternate Hamiltonian algebras. The type-1
Kaplansky algebra is recognized as the derived of the
nonalternate version of the Hamiltonian Lie algebra, the one that
preserves a tensorial 2-form, not an exterior one.

Deforms corresponding to nontrivial cohomology classes can be
isomorphic to the initial algebra, e.g., we confirm Grishkov's
implicit claim and explicitly describe the Jurman algebra as such a
``semitrivial" deform of the derived of the alternate Hamiltonian
Lie algebra. This paper helps to sharpen the formulation of a
conjecture describing all simple finite-dimensional Lie algebras
over any algebraically closed field of nonzero characteristic and
supports a conjecture of Dzhumadildaev and Kostrikin stating that
all simple finite-dimensional modular Lie algebras are either of
``standard" type or deforms thereof.

In characteristic 2, we give sufficient conditions for the known
deformations to be semitrivial.

\end{abstract}

%\date{XXX}

\markboth{\itshape Sofiane Bouarroudj\textup{,} Alexei
Lebedev\textup{,} Dimitry Leites\textup{,} Irina
Shchepochkina}{{\itshape Deformations in characteristic $2$}}

\maketitle

\thispagestyle{empty}

\section{Introduction}

Hereafter, $\Kee$ is an algebraically closed field of characteristic
$p>0$ unless otherwise stated. The letter $p$ also denotes
``momenta" indeterminates but confusion is impossible.

\ssec{Preparatory information} Assuming that $p^\infty=\infty$
and $\Nee=\{1, 2, \dots\}$, we designate
\begin{multline}
\label{u;N} \cO(m; \uN):=\Kee[u;
\uN]:=\Span_{\Kee}\left(u^{(r)}\mid r_i < p^{N_{i}}
\right)\\
\text{~~for $u=(u_1,..., u_m)$ and $r=(r_1,..., r_m)$ and any
$m$-tuple $\uN = (N_1,..., N_m)$},
\end{multline}
where
$N_i\in\Nee\cup\infty$ for any $i$, the addition is natural, and the product is given by
\begin{equation}
\label{divp} u^{(\underline{r})} \cdot u^{(\underline{s})} = \binom
{\underline{r} + \underline{s}} {\underline{r}} u^{(\underline{r} +
\underline{s})}, \text{~~where $\binom {\underline{r} +
\underline{s}} {\underline{r}}:=\prod\limits_{i=1}^m\binom {r_{i} +
s_{i}} {r_{i}}$}.
\end{equation}
The elements of the algebra $\cO(m; \uN)$ of divided
powers serve as ``functions" over $\Kee$. The shearing vector with
smallest coordinates
\begin{equation}\label{n_s}
\uN_s=(1,\dots,1)
\end{equation}
is of particular interest (see \cite{Vi1, BLLS}). Only one of
the algebras of divided powers $\cO(n;\uN)$ is indeed
generated by the indeterminates declared: if $\uN=\uN_s$. Otherwise,
the list of generators consists of $u_i^{(p^{k_i})}$ for all $i$ and
$k_i$ such that $1\leq k_i<N_i$. We define \textit{distinguished}
partial derivatives by setting
\[
\del_i(u_j^{(k)})=\delta_{ij}u_j^{(k-1)}\quad\text{for any }k<p^{N_j}.
\]

Let $\fvect (m; \uN):=\fder_{dist}(\cO(m; \uN))$
be the general \textit{vectorial Lie algebra} spanned by all
distinguished derivations $f_i\del_i$, where $f_i\in\cO(m;
\uN)$; let $\fsvect (m; \uN)$ be its subalgebra
of di\-ver\-gence-free derivations. Various vectorial Lie algebras
are complete or partial \textit{Cartan prolongs}, i.e., the results
of generalized prolongation procedures.

\sssec{Complete Cartan prolongations}\label{CTS} For details, see
\cite{Shch}. Let $DS^{k}$ be the operation of raising to the $k$th
divided symmetric power and
$DS^{\bcdot}:=\mathop{\oplus}\limits_{k\geq 0}DS^{k}$; we set
\begin{equation}
\label{2.5.1'} \begin{array}{l} i\colon DS^{k+1}(\fg_{-1})^*\otimes
\fg_{-1}\tto
 DS^{k}(\fg_{-1})^*\otimes \fg_{-1}^*\otimes\fg_{-1},\\
 j\colon DS^{k}(\fg_{-1})^*\otimes \fg_{0}\tto
DS^{k}(\fg_{-1})^*\otimes \fg_{-1}^*\otimes\fg_{-1}\end{array}
\end{equation}
as the natural maps. Let the \textit{$(k,\uN)$th prolong of the
pair $(\fg_{-1}, \fg_0)$} be
\begin{equation}
\label{genprol} \fg_{k, \uN} = \left
(j(DS^{\bcdot}(\fg_{-1})^*\otimes \fg_0)\cap
i(DS^{\bcdot}(\fg_{-1})^*\otimes \fg_{-1})\right )_{k, \uN},
\end{equation}
where the subscript $k$ in the right-hand side singles out the
component of degree $k$. It is easy to show that if the $\fg_0$-module $\fg_{-1}$ is
faithful, then
$(\fg_{-1},\fg_0)_{*, \uN}=\mathop{\oplus}\limits_k\fg_{k, \uN}$ is
a Lie subalgebra in $\fvect(\dim \fg_{-1}; \uN)$; it is called the (complete)
\textit{Cartan prolong} of the pair $(\fg_{-1},\fg_0)$. A
\textit{partial prolong} is a subalgebra of $(\fg_{-1},\fg_0)_{*,\uN}$
generated by $\fg_{-1}$, $\fg_0$, and a $\fg_0$-submodule of $\fg_1$.

\sssec{Lie algebras of Hamiltonian series}\label{HamAndO} A detailed
description of several types of Hamiltonian series, their
divergence-free subalgebras, their central extensions---Poisson
algebras, and their simple derived in characteristic 2 can be found
in \cite{LeP}. Here, we briefly recall that a given symmetric
bilinear form $B$ on the space $V$ is said to be \textit{alternate}
if $B(v,v)=0$ for any $v\in V$ and \textit{nonalternate} otherwise.
The normal shapes of nondegenerate bilinear forms $B$ whose Gram
matrices are also denoted by $B$ are denoted by $\Pi(n)$ and $I(n)$
 if reduced to the side and main diagonal, respectively, where
$n=\dim V$. The orthogonal Lie algebras $\fo_B(V)$ that preserve
these normal forms are denoted by $\fo_\Pi(V)$ and $\fo_I(V)$,
respectively. If $n$ is odd, then there is only one equivalence
class of nondegenerate symmetric bilinear forms, and we can hence
drop the subscript $B$ in $\fo_B(V)$; if $\dim V$ is even, then
there are two equivalence classes of nondegenerate symmetric forms.

The Hamiltonian Lie algebra can be
\textit{alternate} $\fh_\Pi(V;\uN)$ or \textit{nonalternate}
$\fh_I(V;\uN)$ depending on the type of the differential
2-form\footnote{It was shown in \cite{LeP} that this 2-form is the
sum of products of 1-forms, the product being either exterior
(alternate case) or tensor (nonalternate case).} the algebra preserves
by means of the Lie derivative. We often write $\fh_B(n;\uN)$
instead of $\fh_B(V;\uN)$. Both $\fh_\Pi(V;\uN)$ and $\fh_I(V;\uN)$ have
divergence-free subalgebras described, together with a history of
earlier partial discoveries, in \cite{LeP}.

In this paper, we describe $\fh_B(V;\uN)$ as the Cartan prolong
$(V,\fo_B(V))_{*,\uN}$ with the multiplication given by the
\textit{Poisson bracket} of generating functions
\begin{equation}\label{P.B.}
\{F,G\}_{B}=\mathop{\sum}\limits_{1\leq i,j\leq n}
B_{ij}\pderf{F}{x_i}\pderf{G}{x_j}\quad\text{for any }
F,G\in\cO(n;\uN),\text{ where }(B_{ij})=B.
\end{equation}
The elements of $\fh_B(V;\uN)$ can be realized by vector fields
\begin{equation}\label{ham}
H_{F}=\mathop{\sum}\limits_{1\leq i,j\leq n}
B_{ij}\pderf{F}{x_i}\pder{x_j}\quad\text{for any } F \in\cO(n;\uN),
\end{equation}
where $(B_{ij})_{i,j=1}^n=B$. The Lie algebra whose space is
$\cO(n;\uN)$ with the bracket \eqref{P.B.} is called the
\textit{Poisson algebra} if $B\not\sim I$; it is a central extension
of $\fh_B(n;\uN)$ for $B\sim\Pi$.

We note that although $\fh_I(n;\uN)$ is well defined, \textbf{there is
no Lie algebra} $\fpo_I(n;\uN)$ with the bracket \eqref{P.B.}.
Indeed, the bracket should be antisymmetric, i.e., alternate,
while $\{x_i,x_i\}_{I}=1$, not 0. More on possible brackets
corresponding to the alternate bilinear form $B$ can be found in
\cite{LeP} and subsec.~\ref{Pb}.

\ssec{Overview of the situation} Even the incomplete stock of
nonisomorphic spe\-cies in the zoo of simple finite-dimensional Lie
algebras for $p=2$ was until recently considered uncomfortably
numerous (see the introduction in \cite{S}). It has many more exhibits
than would have been considered ``normal" if the classification in
cases $p>3$ were taken as a ``norm." A final touch in the
proof of the classification can be found in \cite{BGP}.

The improved version of the Kostrikin--Shafarevich conjecture due to
Dzumadildaev and Kostrikin \cite{KD} states that for any $p>0$,
\begin{equation}\label{KDconj}
\begin{minipage}[l]{12cm}
\textbf{all simple Lie algebras are either of ``standard" type\\ or
deforms (the results of deformations) thereof}.
%\end{center}
\end{minipage}
\end{equation}
The improved conjecture definitely embraces $p\geq 5$ as proved in
\cite{KD}, where the Melikyan algebras are identified as deforms of
the Poisson algebras. The claim in \cite{KuJa} identifying the
simple Ermolaev algebras as deforms of the contact algebras supports
the conjecture \eqref{KDconj} for $p=3$. The conjecture
\eqref{KDconj} seems plausible if, for $p=2$ and 3, we enlarge the
stock of examples ``standard" for $p\geq 5$ by exhibits from
\cite{GL, SkT1, BGL1, BGL3, LeP, BGLLS,BGLLS1}, mostly found after
\eqref{KDconj} was formulated. Therefore, we must find out which
simple Lie algebras are ``standard" from the standpoint of
\eqref{KDconj}, and solve a ``small technical problem" of describing
all nonisomorphic deforms.

A conjectural description of ``standard", hence, all (together with
their deforms) simple finite-dimensional Lie algebras over $\Kee$
for $p=2$, although longer than that for $p>3$, is possible to
grasp. This conjecture stemmed from an idea that had already led to
the classification of simple Lie superalgebras of polynomial vector
fields over $\Cee$ (see \cite{LSh1}). The new conjecture yielded new
examples for $p=3$ (see the arXiv version of \cite{GL}).

For $p=2$, the new conjecture (see \cite{Ltow} for a briefly
formulated version) gathers all examples known to us in describable
groups and indicates ways to obtain new examples. Apart from these
and several new examples given in \cite{BGLLS1}, there are also
known examples (due to Kaplansky, Shen, Skryabin, Brown, Jurman,
Vaughan-Lee, and recently Eick) of mysterious nature. In this paper
and \cite{BGLLS, BGLLS1}, we study which of these ``mysterious"
examples, if any, might qualify as ``standard" from the standpoint
of Conjecture \eqref{KDconj}, see also \cite{SkT1, GZ}. We demystify
the other examples by identifying them either as deforms of or as
isomorphic to some of the ``standard" examples.

\sssec{On limited information derived from cohomology in describing
deforms of Lie algebras}\label{semitr} In \S\ref{BLW}, we recall how
deformations of Lie algebras are calculated. The \textit{trivial
deformation} of $\fg$ corresponds to the change of the basis in
$\fg$ corresponding to a 2-coboundary, while the linear part of any
global deformation is a cocycle; hence, deforms linear in the
deformation parameter, also called \textit{infinitesimal} deforms,
correspond to cocycles representing classes of $H^2(\fg;\fg)$. If
$\Char\Kee>0$, then there are ``fake deformations," which means not
that some of linear deforms corresponding to cocycles representing
classes of $H^2(\fg;\fg)$ might be not extendable to a global
deformation but something much worse. The textbooks and papers on
Lie (super)algebra cohomology do not yet indicate the following
important phenomenon:
\begin{equation}\label{fake}
\begin{minipage}[l]{12cm}
\textbf{Let each cocycle representing a class of $H^2(\fg;\fg)\neq0$
be extendable to a global deformation. This does not preclude some
(or all) deforms of $\fg$ from being isomorphic to $\fg$}.
%\end{center}
\end{minipage}
\end{equation}
Let \textit{semitrivial} deformations (and their results, the
deforms) be the ones whose linear parts are given by cocycles
representing nontrivial cohomology classes but whose deforms are
isomorphic to the initial Lie algebra. In addition to examples in
\cite{BLW}, we show that the Jurman algebras are \textit{semitrivial
deforms}.

Examples of semitrivial deforms have been known to us since 1987 when
we computed that\footnote{Hereafter, the prime ${}'$ denotes the first
derived, also called the commutant.} $\dim H^2(\fo'(3); \fo'(3))=2$ over
$\Kee$ for $p=2$ while, up to an isomorphism, there is only one
simple 3-dimensional Lie algebra: $\fo'(3)$. (Ten years earlier the
phenomenon \eqref{fake} was observed without any explanations of its
origin in \cite{DzhK}.) The first \textit{explanation} of the cause
of the phenomenon \eqref{fake} was given in \cite{BLW}. We describe
a certain type of semitrivial deformations for $p=2$ in
subsec.~\ref{semit}.

\sssec{The Vaughan-Lee algebras are not new over $\Kee$} The table
on p.~948 in \cite{Ei} shows that simple algebras of Vaughan-Lee
(all new over $\Fee_2$) are only new as forms of Lie algebras known
over $\Kee$ (or even over a Galois field extending $\Fee_2$).

\sssec{The Eick algebras are new} Several simple Lie algebras were
introduced in \cite{Ei} and conjectured (e.g., because the list of
``known" algebras Eick used for comparison was incomplete as
compared with a wider list known to us) to be new. These algebras
had to be interpreted and described in more detail than in
\cite{Ei}. With Eick's help, we recently established that all the
six tentatively new algebras in \cite{Ei} are indeed new. All the
six new Eick algebras are obtained in one of the ways predicted by
the conjecture \cite{Ltow}: Eick algebras are \textit{partial}
Cartan prolongs (see subsec.~\ref{CTS}), like Frank algebras for
$p=3$ (cf.~\cite{GL}), or deforms of something ``standard."

\sssec{One of the Shen algebras and its generalized Cartan prolong
due to Brown are ``standard"} In \cite{Sh}, Shen described several
simple Lie algebras. One Shen algebra was rediscovered, together
with several algebras new at that time, by Brown \cite{Bro}. Brown's
examples, described only in components in \cite{Bro}, were
interpreted in \cite{GL, BGLLS} together with a clarification of
their structure and related new simple Lie superalgebras. One
remarkable exceptional simple Lie algebra (we call it $\fg\fs(2)$ in
honor of Shen Guangyu, who discovered it; see \cite{Sh}) is a true
analog of the Lie algebra $\fg(2)$ in characteristic 2. Brown
rediscovered this algebra (Eick called it $\text{Bro}_2(1,1)$ in
\cite{Ei}) and considered its Cartan prolong. We called the derived
of this prolong $\fme'(5;\uN)$, and interpreted it as a version of
the Melikyan algebras in characteristic 2 (see \cite{BGLLS1}). These
$\fme'(5;\uN)$, where $N_1=N_2=N_3=1$ are the only possible values,
and its particular case $\fg\fs(2)=\fme'(5;\uN_s)$ seem to be new
``standard" examples. Several Shen algebras were interpreted in
\cite{BGL2} as deforms of certain ``standard" algebras, and several
of Shen's examples are either nonsimple or not new (see \cite{LLg}).
Moreover, the multiplication in several Shen algebras does not
satisfy the Jacobi identity, and we could not repair this.

\sssec{Jurman and Kaplansky algebras as deforms} We started this
paper intending to prove that the Jurman and Kaplansky algebras are
deforms of more ``conventional" simple Lie algebras (such as the
two nonisomorphic versions of the Lie algebra of Hamiltonian vector
fields) and their divergence-free subalgebras (see \cite{LeP}). While
this paper was being written, Grishkov published a note\footnote{A draft
of this note was available on
Grishkov's home page since 2009.} \cite{GJu} claiming that the
Jurman algebra is \textbf{isomorphic} to the (derived of) a Hamiltonian
Lie algebra. Grishkov's paper is based on a difficult result due to
Skryabin, and its claim is implicit. Consequently, we heard doubts
that it is correct. It IS correct: we give an explicit isomorphism
in Prop.~\ref{Jcocycles}. Amazingly, the existence of this
isomorphism does not contradict the fact that the Jurman algebra is
a deform corresponding to a cocycle personifying a nontrivial
cohomology class of the (derived of the) Hamiltonian Lie algebra:
Jurman algebras are examples of ``semitrivial" deforms (see
subsec.~\ref{semitr} and~\ref{semit}.

In \S\ref{secK}, we identify type-1 Kaplansky algebras with certain
known ``standard" Lie algebras and prove that type-2 Kaplansky
algebras are deforms of certain ``standard" Lie algebras. Type-3
algebras were identified (in different terms) by Kaplansky himself
as $\fo'_I(n)$. Type-4 algebras might be ``standard", see Subsection
\ref{KapNotO}.

\ssec{Main results}\label{ssmainres} The three \textbf{most
interesting parts} of our paper:

(1) The discovery of a $\Zee/2$-grading quadratic in roots. Among
the Lie algebras known to us, the type-2 and type-4 Kaplansky
algebras are the only ones with such gradings (cf.~\eqref{KapS}).
These Kaplansky algebras are unique among the Lie algebras known to
us. We present details on relations between gradings and
derivations, in particular, an observation that gradings are not
always defined by derivations, in Subsection~\ref{grader}.

(2) A method assigning a simple Lie superalgebra to \textbf{every}
simple Lie algebra. It is illustrated with seven new series,
superizations of the Kaplansky algebras in
Subsection~\ref{SSgenSuper}.

(3) Thanks to the insistence of a referee, we managed to give
sufficient conditions for deformations, encountered so far for
$p=2$, to be semitrivial (Theorem~\ref{semiTr}); but other types of
semitrivial deformations are possible.

The \textbf{main bulk} of the paper is devoted to interpreting the
simple Lie algebras discovered by Jurman and Kaplansky in terms of
better known (``standard") examples of Lie algebras of Hamiltonian
vector fields or their simple derived. Voluminous computations are
performed using Grozman's \textit{Mathematica}-based package
\textit{SuperLie}, cf. Lemmas \ref{conjlemma}, \ref{Lbn},
\ref{Lbn2}, \ref{nonlin}.

\ssec{Open problems}\label{OP}
\textbf{1.} We have described deforms of $\fh_\Pi'(n;\uN)$ for $n=2$.
Investigation
of the isomorphism classes of the deforms for $n>2$ and any $\uN$ is
a must. The classification of the deforms of the more natural
nonsimple relative of the simple algebra, i.e., of $\fh_\Pi(n;\uN)$,
is also needed: it had led to an interpretation of previously
mysterious type-2 Kaplansky algebras. The search for deforms of
another relative, the Poisson Lie algebra $\fpo_\Pi(n;\uN)$, is
equally reasonable (answers to such problems for $p=0$ have
physical interpretations), see \cite{KT, KD}.

\textbf{2.} A new way to construct simple Lie algebras in the
absence of classification\footnote{Eick's search
is random; to estimate the
probability of a miss is very interesting.} is provided
in \cite{Ei}  if the structure constants belong to
$\Fee_2$; the parametric families cannot be captured by Eick's
method. Although regrettably restricted to algebras of small
dimension (currently $\leq 20$), Eick's computer-aided approach
(when its range will have been widened to dimension 250,
or at least 80) promises to give a base for a
conjecture making its theoretical proof psychologically comfortable.

\textbf{3.} Our results show that $\fh_\Pi(2;\uN)$ and
$\fh_\Pi'(2;\uN)$ have different numbers of deforms and both types
of deforms are important for classifying simple Lie algebras. The
situation is similar to that in characteristic 0, where the Lie
superalgebra $\fh(2n|m)$ has (in the only case $(2n|m)=(2|2)$) more
deformations than $\fpo(2n|m)$ (see \cite{LSh2}). A problem is to
describe deformations of $\fh_\Pi(2n;\uN)$ and $\fh_\Pi'(2n;\uN)$.

\textbf{4.} In \cite{Sk}, Skryabin classified nonequivalent normal
shapes of the exterior 2-forms preserved by the Hamiltonian Lie algebras
$\fh_\Pi(n;\uN)$. It remains to classify \`{a} la \cite{Sk} nonequivalent
normal shapes of the tensorial 2-forms preserved by the Hamiltonian Lie
algebras $\fh_I(n;\uN)$ and its divergence-free subalgebras.

\textbf{5.} Listing all nonisomorphic deforms of
$\fg=\fh_\Pi'(2;(g,h+1))$ requires considering orbits of the
$\Aut(\fg)$-action on the space $H^2(\fg;\fg)$ following Kuznetsov
and his students (see \cite{KCh,Ch}). If the algebraic group
$\Aut(\fg)$ of automorphisms of $\fg$ were computed, then the result
in \cite{FG} could be extended to the simple Lie algebras without a
Cartan matrix. So far, this has been done only in certain particular
cases (see \cite{Pre}).

\textbf{6.} Several identification problems remain (see
subsecs.~\ref{conjlemma}, \ref{parKapS}, \ref{Prob}, \ref{sssIso}, and
\ref{gh=31prob} and eq.~\eqref{minip}).

\textbf{7.} A tough problem, partly solved in \S\ref{semit}, is to
to give sufficient conditions for a  deformation to be semitrivial.

\section{Deformations and cohomology (\cite{BLW})}\label{BLW}

\subsection{Lie algebras} A \textit{multiparameter deformation} of a Lie algebra $\mathfrak{g}$
over $\Kee$ is a Lie algebra $\mathfrak{g}_{t}$, where
$t=(t_1,\ldots,t_r)$, given by a Lie algebra structure on the tensor
product $\mathfrak{g}\otimes_\Kee\Kee[[t]]$ such that the Lie
algebra $\mathfrak{g}_{0}$, i.e., the algebra obtained for $t=0$, is
isomorphic to $\mathfrak{g}$. For any $x,y\in\mathfrak{g}$, the deformed bracket has the form
\begin{equation}\label{1vykl}
\begin{split}
 [x,y]_{t_1,\ldots,t_r}={}&c^0(x,y)+t_1c_1^1(x,y)+\ldots+t_r c_r^1(x,y)\\
 &{}+t_1^2c^2_{1,1}(x,y)+t_1t_2c^2_{1,2}(x,y)\ldots+t_r^2
c^2_{r,r}(x,y) +\ldots ,
\end{split}
\end{equation}
 where $c^0(x,y):=[x,y]$. By linearity, it
suffices to specify the deformed bracket of elements in $\mathfrak{g}$.
The degree-1 conditions say that the maps
$c_i^1\colon\mathfrak{g}\otimes_\Kee\mathfrak{g}\tto\mathfrak{g}$
must be antisymmetric and must be $2$-cocycles, i.e., for all
$i=1,\ldots,r$ and any $x,y,z\in\mathfrak{g}$, we have
\[
\begin{split}
 dc_i^1(x,y,z):={}&c_i^1([x,y],z)+c_i^1([y,z],x)+c_i^1([z,x],y)\\
 &{}-[x,c_i^1(y,z)]-[y,c_i^1(z,x)]-[z,c_i^1(x,y)]=0.
\end{split}
\]
For brevity, we recall properties of deformations for $1$-parameter
deformations; the multidimensional case is considered routinely.
For example, eq.~\eqref{1vykl} becomes
\begin{equation}\label{2vykl}
 [x,y]_{t}=c^0(x,y)+tc^1(x,y)+t^2c^2(x,y)+\ldots.
\end{equation}
Two (formal) $1$-parameter deforms $\mathfrak{g}_t$ and
$\tilde{\mathfrak{g}}_t$ given by the collections $c=(c^1, c^2,
\dots)$ and $\tilde{c}=(\tilde{c}^1, \tilde{c}^2, \dots)$ lead to
equivalent deforms (i.e., $\mathfrak{g}_t$ and
$\tilde{\mathfrak{g}}_t$ are isomorphic as Lie algebras by an
isomorphism of the form $\tau(x;t)=\mathop{\sum}\limits_{i\geq 0}\tau_i(x)t^i$,
where $\tau_0=\id$, for any $x\in\mathfrak{g}$) if and only if $c$
and $\tilde{c}$ are related as follows (for all $n>0$, where $i,j,k\geq 0$):
\begin{equation}\label{trivDef}
 \mathop{\sum}\limits_{i+j=n}\tau_i(\tilde{c}^j(x,y))=
 \mathop{\sum}\limits_{i+j+k=n}c^i(\tau_j(x), \tau_k(y)).
\end{equation}
Any change of basis of $\fg$ can be included in a 1-parameter family
$\tau(\cdot;t):\fg\tto\fg$ and regarded, naturally, as a trivial
deformation. We see that the trivial deformation corresponds to
$d\tau_1$ modulo $t^2$; the search for the most general
multiparameter deformation of a given Lie algebra therefore begins
with computing the space $H^2(\mathfrak{g};\mathfrak{g})$. Its
explicit basis given by $2$-cocycles (representing the classes)
determines infinitesimal deformations. We then try to prolong each
infinitesimal deformation to higher degrees. The Jacobi identity
imposes conditions on all terms in the deformed bracket, which must
be satisfied in each degree. In particular, two $1$-parameter
degree-1 cocycles $c^1$ and $\tilde{c}^1$ are
\textit{infinitesimally equivalent} (i.e., $\tau={\rm id}+t\tau_1$
modulo $t^2$) if and only if $c^1-\tilde{c}^1=d\tau_1$.

Let $\mathfrak{g}_t$ be a $1$-parameter deformation of a Lie algebra
$\mathfrak{g}$ given by the collection ${c=(c^1, c^2,\dots)}$. The
Jacobi identity yields that the coefficient of $t^{n}$ vanishes for
each $n$:
\begin{equation} \label{br}
 \mathop{\sum}\limits_{0\leq i,j\leq n;\;\; i+j=n}(c^i(c^j(x,y),z)+{\rm
 cyclic}(x,y,z))=0,
\end{equation}
where ${\rm cyclic}(x,y,z)$ denotes the sum of all cyclic permutations
of the arguments of the expression written on the left of it. We set
\begin{equation}\label{Massey}
\begin{array}{l}
[[c^i,c^j]](x,y,z):=c^i(c^j(x,y),z)+c^j(c^i(x,y),z)+{\rm  cyclic}(x,y,z),\\
c^k\circ c^k (x,y,z):=c^k(c^k(x,y),z)+{\rm  cyclic}(x,y,z).
\end{array}\end{equation}
The brackets $[[c^i,c^j]]$ are called \textit{Nijenhuis brackets} (in
differential geometry) or \textit{Massey brackets} (in deformation
theory). The sum \eqref{br} can be expressed as a
\textit{Maurer--Cartan} equation:
\begin{equation}\label{MaurCart}
dc^n=\sum_{0< i<j\leq n;\;\; i+j=n}[[c^i,c^j]]+
\begin{cases}0&\text{if }n=2k+1,\\
 c^k\circ c^k&\text{if }n=2k.\end{cases}
\end{equation}

To prolong an infinitesimal deformation given by a cocycle $c^1$, we
first compute $[[c^1,c^1]]$. If $[[c^1,c^1]]=0$, then the infinitesimal
deformation satisfies the Jacobi identity and is a true deformation.
If $[[c^1,c^1]]\in Z^3(\mathfrak{g},\mathfrak{g})$ and
$[[c^1,c^1]]\not\in B^3(\mathfrak{g},\mathfrak{g})$, then the
infinitesimal deformation is obstructed and cannot be prolonged.
If $[[c^1,c^1]]=dc^2$ with $c^2\not=0$, then $-t^2c^2$ is the
second-degree term of the deformation. To prolong the deformation to
degree 3, we compute the Massey product $[[c^1,c^2]]$. There are the
three possibilities
\[
1)\;[[c^1,c^2]]=0,\quad
2)\;[[c^1,c^2]]=dc^3\text{ for some }c^3\not=0,\quad
3)\;[[c^1,c^2]]\not=dc^3\text{ for any }c^3.
\]
If $[[c^1,c^2]]=dc^3$, then $-t^3c^3$ gives the third-degree
prolongation of the deformation. To go up to degree $4$, we must be
able to compensate $c^2\circ c^2+[[c^1,c^3]]$ by a coboundary $dc^4$,
and so on. The main difficulty here is that the representatives of
the cohomology classes and the cochains  $c^2$, $c^3$, etc., are not
uniquely\footnote{If $c^2$ is a solution of the equation
$dc^2=c^1\circ c^1$, then $c^2+\mathrm{cocycle}$ is also a solution.
The choice of a certain $c^2$ affects the expression of the $c^3$
terms. The problem is how to find a ``nice" $c^2$ in order to have
as few $c^3$ terms as possible and, more importantly, vanishing
Massey products in degrees $>3$. If we fail to achieve this with
$c^2$, then we can try with $c^3$, and so on.} defined. A good
choice of cochains may considerably facilitate computations.
The following lemma is helpful.

\sssec{Grozman's lemma (\cite{BLW})}\label{LGr}
For any finite-dimensional Lie algebra $\mathfrak{g}$, all cochains
with coefficients in the adjoint module can be expressed as sums of
tensor products of the form $a\otimes\omega$, where $a\in\mathfrak{g}$
and $\omega\in\bigwedge^{\bcdot}(\mathfrak{g}^*)$.

\begin{Lemma}\label{lem1} For any $c=a\otimes\omega$ with
$a\in\mathfrak{g}$ and $\omega\in\bigwedge^{r}(\mathfrak{g}^*)$,
let $dc$ denote the coboundary of $c$ in the complex with
coefficients in the adjoint module, $d\omega$ be the coboundary in
the complex with coefficients in the trivial module, and $da$ be the
coboundary of $a\in\mathfrak{g}$ regarded as a $0$-cochain in the
complex with coefficients in the adjoint module. Then
$dc=a\otimes d\omega + da\wedge\omega$.
\end{Lemma}

\subsection{Lie superalgebras} In certain problems, if $p>0$, then we must replace the
(formulas of) conventional cohomology in the preceding subsection
with divided power ones, see \cite{BGLL}.

\subsection{Semitrivial deformations}\label{semit} In all examples we
know, the deforms of $\fg$ corresponding to semitrivial deformations
for $p=\Char\Kee$ are isomorphic to the initial Lie algebra via an
isomorphism $\tau(\cdot,t)$ given by an expression of the form (also
see \cite{BLW}, where equality $\sqrt[p]{(1+t)^i}=1+(\sqrt[p]{t})^i$ should
be used)
\begin{equation}\label{semiteq}
\tau(x;t)=x+\mathop{\sum}\limits_{i\geq 1}
\tau_i(x)(\sqrt[p]{t})^i\quad\text{with the }\tau_i
\text{ satisfying conditions \eqref{trivDef}}.
\end{equation}
Therefore, the semitrivial cocycle $c$, that would have been
obtained for $p=0$ as the differential of a $1$-cochain $C$ such
that $C(x)=\pderf{\tau(x,t)}{t}$ for any $x\in\fg$, cannot be thus
obtained if $p>0$, because the function $\sqrt[p]{t}$ is not
differentiable for $p>0$.

For any derivation $D$ of the Lie algebra $\fg$, let
$c_D \in C^2(\fg;\fg)$ be defined by
\begin{equation}\label{c_D}
c_D(x,y) = [Dx,Dy]\quad\text{for any }x,y\in\fg.
\end{equation}
It is easy to verify that $c_D$ is a cocycle, i.e., $dc_D(x,y,z) = 0$
for all $x,y,z\in\fg$. Moreover, if $p\neq 2$, then
$c_D=-\frac12d(D^2)$, i.e., $c_D$ is trivial.

Here is a general description of semitrivial deformations for $p=2$
in the case where an isomorphism between the algebra and its deform
is a polynomial in $\sqrt{t}$. We do not know how to characterize
arbitrary semitrivial deformations.

\sssbegin{Lemma}\label{semiTr0} Let $\fg$ be a Lie algebra in
characteristic $2$, and let $c\in Z^2(\fg;\fg)$ generate a semitrivial
deformation such that the isomorphism between the algebra and
its deform is polynomial in $\sqrt{t}$. Then there is a derivation
$D$ of $\fg$ such that $c$ is equivalent to $c_D$, i.e.,
$c-c_D\in B^2(\fg;\fg)$.
\end{Lemma}

\begin{proof} Let the isomorphism be $F_t(x)=x+\sum\limits_{i\geq 1}
f_i(x)(\sqrt{t})^i$, where only a finite number of the $f_i\in
C^1(\fg;\fg)$ are nonzero. Then the deformed bracket has the form
\begin{equation}\label{defBr}
\begin{split}
[x,y]_t ={}& [x,y] + (f_1([x,y])+[f_1(x),y]+[x,f_1(y)])\sqrt{t}
 \\
&{}+([f_1(x),
f_1(y)] + df_2(x,y))t+\text{higher order terms in }\sqrt{t}.
\end{split}
\end{equation}
Because the coefficient of $\sqrt{t}$ is zero and the coefficient
of $t$ is $c(x,y)$, we see that $f_1$ is a derivation, and $c$ is
equivalent to $c_{f_1}$.
\end{proof}

We do not know whether every cocycle of the form (\ref{c_D})
generates a global deformation, nor do we know whether all such
deformations are semitrivial or trivial. We can prove it for two
classes of derivations.

\sssbegin{Theorem}\label{semiTr} Let $\fg$ be a Lie algebra over a
perfect field $\Kee$ of characteristic $2$, and let $D$ be a derivation
of $\fg$ satisfying at least one of the following two conditions:

1. $D^2=\alpha D$ for some $\alpha\in\Kee$ (if $\alpha=1$, then such
$D$ corresponds to a $\Zee/2$-grading);

2. $D$ is nilpotent.

Then the cocycle $c_D$ (see Lemma $\ref{semiTr0}$) is either trivial
or semitrivial.
\end{Theorem}

\begin{proof} To prove that $c$ is either trivial or semitrivial, we
show that there is a polynomial family of invertible maps
$F_t(x) = x + \sum\limits_{i\geq 1}f_i(x)t^i$, where only a finite
number of the $f_i\in C^1(\fg;\fg)$ are nonzero, such that
\begin{equation}\label{stdef}
{}[F_t(x), F_t(y)] = F_t([x,y] + \sum_{j\geq 1} c_j(x,y)t^{2j})\quad
\text{for all }x,y\in\fg,
\end{equation}
where only a finite number of the $c_j\in C^2(\fg;\fg)$ are
nonzero, and $c_1 = c_D$. This would mean that there exists a
global deformation of $\fg$ given by the bracket
\[
{}[x,y]_t = [x,y] + \sum_{j\geq 1} c_j(x,y)t^{j} ,
\]
all the deformed algebras being isomorphic to $\fg$ with the
isomorphism given by $F_{\sqrt{t}}$.

1. If $D^2=\alpha D$, then we set $F_t(x) = x + tD(x)$; this map is
invertible if $t\alpha\neq 1$. We also set $c_1 = c_D$ and $c_j = 0$
for $j>1$. It is easy to see that (\ref{stdef}) holds because
$D[Dx,Dy]=0$.

2. Now let $D$ be nilpotent. We set
\[
f_i = \begin{cases} D^i & \text{if }i=2^k\text{ for some }k,\\
0 &\text{otherwise.}\end{cases}
\]
We note that all the $f_i$ are derivations. We set\footnote{It might be
wondered, ``Why does the deformed bracket satisfy the Jacobi
identity? In particular, why do the $c_j$ satisfy the Maurer--Cartan
equations?" We forget for a moment that $\fg$ is a Lie algebra and
regard it just as an algebra. Then if \eqref{stdef} holds (and we show
that it does), then the deformed algebra is isomorphic to $\fg$ with
the isomorphism given by $F$. Because $\fg$ is a Lie algebra, the
deformed algebra is also a Lie algebra, i.e., the deformed bracket
satisfies Jacobi identity.}
\[
c_j(x,y)=\begin{cases}\sum\limits_{1\leq s\leq 2j-1}[D^sx,D^{2j-s}y]
& \text{if }j=2^k\text{ for some }k,\\
0 & \text{otherwise.}\end{cases}
\]

Because $D$ is nilpotent, the map $F_t$ is invertible. We must prove
that the coefficients of $t^n$ in \eqref{stdef} are the same for any
$n>0$, i.e., that
\begin{multline}
\label{new1}
{}[f_n(x),y] + [x, f_n(y)] + \sum\limits_{1\leq i\leq n-1}[f_i(x),
f_{n-i}(y)]\\
=f_n([x,y]) +\begin{cases}
c_m(x,y)+\sum\limits_{1\leq j\leq m-1} f_{n-2j}(c_j(x,y))&\text{if }n=2m,\\
\sum\limits_{1\leq j\leq m} f_{n-2j}(c_j(x,y))&\text{if }n=2m+1.
\end{cases}
\end{multline}
Because $f_n$ is a derivation, eq.~\eqref{new1} can be simplified to
\begin{equation}\label{stdef2}
\sum\limits_{1\leq i\leq n-1}[f_i(x), f_{n-i}(y)] = \begin{cases}
c_m(x,y)+\sum\limits_{1\leq j\leq m-1} f_{n-2j}(c_j(x,y))
 &\text{if }n=2m,\\
\sum\limits_{1\leq j\leq m}f_{n-2j}(c_j(x,y)) &\text{if }n=2m+1.
\end{cases}
\end{equation}

If $n=1$ or if $n$ cannot be represented as a sum of at most two powers
of $2$, then both sides of \eqref{stdef2} are equal to $0$ because $f_i$
and $c_i$ are only nonzero if $i$ is a power of 2. Therefore,
\eqref{stdef2} has nonzero terms only if $n$ is the sum of at most
two powers of 2.

If $n = 2^k+2^l$, where $k>l\geq 1$, then the right-hand side of
\eqref{stdef2} is
\begin{multline*}
f_{2^k}(c_{2^{l-1}}(x,y)) + f_{2^l}(c_{2^{k-1}}(x,y))\\
\begin{aligned}
={}&\sum\limits_{1\leq s\leq 2^l-1}\!([D^{2^k+s}x, D^{2^l-s}y] +
 [D^sx, D^{2^k+2^l-s}y])\\
&{}+\sum\limits_{1\leq s\leq2^k-1}\!([D^{2^l+s}x,
D^{2^k-s}y] + [D^sx, D^{2^k+2^l-s}y])\end{aligned}\\
=[D^{2^l}x, D^{2^k}y] + [D^{2^k}x, D^{2^l}y],
\end{multline*}
and the last expression is exactly the left-hand side of \eqref{stdef2}.

If $n=2^k+1$, where $k\geq1$ (this case is in fact a particular subcase
of the above case for $l=1$, but we consider it separately because there
is no $c_{2^{l-1}}$ in this case), then the right-hand side
of \eqref{stdef2} is
\[
\begin{aligned}
f_1(c_{2^{k-1}}(x,y)) &=\sum\limits_{1\leq s\leq 2^k-1}([D^{s+1}x,
D^{2^k-s}y] + [D^sx, D^{2^k+1-s}y])\\
&= [Dx, D^{2^k}y]+[D^{2^k}x,Dy],
\end{aligned}
\]
and the last expression is exactly the left-hand side of \eqref{stdef2}.

If $n = 2^k$, where $k\geq 2$, then the right-hand side of \eqref{stdef2}
is
%\[
\begin{multline*}
f_{2^{k-1}}(c_{2^{k-2}}(x,y))+%{}&
c_{2^{k-1}}(x,y)\\
=%{}&
\sum\limits_{1\leq s\leq 2^{k-1}-1}([D^{2^{k-1}+s}x, D^{2^{k-1}-s}y]+
[D^sx, D^{2^k-s}y])+\sum\limits_{1\leq s\leq 2^k-1} [D^sx,D^{2^k-s}y]\\
=%{}&
{}[D^{2^{k-1}}x, D^{2^{k-1}}y],
\end{multline*}
%\]
and the last expression is exactly the left-hand side of \eqref{stdef2}.

Finally, if $n=2$, then the left-hand side of \eqref{stdef2} is
$[Dx,Dy]$, same as its right side, $c_1(x,y)$.
\end{proof}

\section{Modular vectorial Lie algebras as deforms of each
other}\label{sdefs}

Weisfeiler and Kac were the first to discover parametric families of
simple finite-dimensional Lie algebras with a Cartan matrix over
$\Kee$ (see \cite{WK}). For further examples of deforms of simple Lie
algebras, see \cite{DzhK,Dzh,Sk,KD,KuJa,GL,BLW,LeP}.

In this section, we extend the list of such examples and also show
that several nonisomorphic Poisson Lie algebras are deforms of one
Lie algebra nonsimple over $\Kee$ but simple over a ring and thus
resemble forms over algebraically nonclosed fields of an algebra
defined over an algebraically closed field. We regard expressions of the form $k\mod p$,
where $k\in\Zee$, as integers from the segment $[0,p-1]$ and not as
elements of $\Kee$.

\ssbegin{Lemma}\label{LemmaPrep} We consider a linear endomorphism
$\Phi_\alpha$, where $\alpha\in\Kee$, of the algebra $\cO(1;
\underline{n})$, given by the formula
\begin{equation}\label{}
\Phi_\alpha(x^{(k)}) = \alpha^{\left[\frac{k}{p}\right]} x^{(k)},
\end{equation}
where $\left[\nfrac{k}{p}\right]$ denotes the integer part of
$\nfrac{k}{p}$ and $k<p^n$. If $\alpha\neq 0$, then $\Phi_\alpha$ is
an automorphism of $\cO(1; \underline{n})$.\end{Lemma}

\begin{proof} Clearly, $\Phi_\alpha$ is a bijection, and we need
only prove that
\begin{equation}\label{}
\Phi_\alpha(x^{(k)}\cdot x^{(l)}) = \Phi_\alpha(x^{(k)})\cdot
\Phi_\alpha(x^{(l)}),
\end{equation}
i.e.,
\begin{equation}\label{}
\alpha^{\left[\frac{k+l}{p}\right]}\binom{k+l}{k} x^{(k+l)} =
\alpha^{\left[\frac{k}{p}\right]+\left[\frac{l}{p}\right]}
\binom{k+l}{k} x^{(k+l)}.
\end{equation}
One can see that\footnote{The thing equal to $0$ in the
second line of \eqref{twolines} is not the same as the thing equal
to $\left[\frac{k}{p}\right]+\left[\frac{l}{p}\right]$ in the first line.
We also note that in the first line, the equality (involving integer parts
 used as power degrees)
is over $\Zee$; in
the second line (involving the binomial coefficient), the equality is
over $\Kee$ or modulo $p$. In both lines, the residues of $k$ and
$l$ modulo $p$ should be understood as integers from the segment
$[0, p-1]$; then the inequalities make sense.}
\begin{equation}\label{twolines}
\begin{alignedat}{2}
&\left[\frac{k+l}{p}\right]=\left[\frac{k}{p}\right]+
\left[\frac{l}{p}\right]&\quad&\text{if }(k\bmod p)+(l\bmod p)<p,\\
&\binom{k+l}{k} \equiv 0(\bmod p)&\quad&\text{if }(k\bmod p)+(l\bmod
p)\geq p,\end{alignedat}
\end{equation}
and the statement in the lemma hence holds in both cases.\end{proof}

We consider the endomorphism $D_\alpha = \Phi_\alpha^{-1}\circ \del
\circ \Phi_\alpha$ of $\cO(1; \underline{n})$ given explicitly by
\begin{equation}\label{D_alpha}
D_\alpha(x^{(k)}) = \begin{cases}
\del x^{(k)} &\text{if }p\nmid k;\\
\alpha\del x^{(k)} &\text{if }p\mid k.\end{cases}
\end{equation}
We define $D_0$ (i.e., $D_\alpha$ for $\alpha=0$, when $\Phi_0$ is
not defined) using relation \eqref{D_alpha}.

We note that if we consider the isomorphism between $\cO(1;
\underline{n})$ and $\cO(2;(1,n-1))$ given by
\begin{equation}\label{}
x^{(k)} \longleftrightarrow y_1^{(k\bmod p)}
y_2^{\left(\left[\frac{k}{p}\right]\right)},
\end{equation}
then $D_0$ on $\cO(1; \underline{n})$ corresponds to $\del_1$ on
$\cO(2;(1,n-1))$.

Similarly, in the algebra $\cO(d;\uN)$ with indeterminates $x=(x_1,
\dots, x_d)$, we can consider the map
\begin{equation}\label{}
\Phi_\alpha(x^{(r)}) = \alpha^{\mathop{\sum}\limits_{1\leq i\leq d}
\left[\frac{r_i}{p}\right]}x^{(r)},
\end{equation}
which is an isomorphism for $\alpha\neq 0$. The maps (here
$\del_i:=\del_{x_i}$)
\begin{equation}\label{D_a}
D_{\alpha,i} = \Phi_\alpha^{-1}\circ \del_i \circ \Phi_\alpha\quad
\text{ act as }D_{\alpha,i}(x^{(r)}) = \begin{cases} \del_i
x^{(r)} &\text{if }p\nmid r_i,\\
\alpha\del_i x^{(r)} &\text{if }p\mid r_i.\end{cases}
\end{equation}
We define $D_{0,i}$ using the relations \eqref{D_a}.

\ssec{Poisson Lie algebras}\label{Pb} We consider the Lie algebra
$\fpo_B(d;\uN)$, where $B=(B_{ij})$ is an alternate (the analog of
antisymmetric for $p=2$) nondegenerate bilinear form on a
$d$-dimensional space. The space of this algebra coincides with
$\cO(d;\uN)$, and the Poisson bracket is defined by eq.~\eqref{P.B.}.

We consider the deformed bracket of $\fpo_B(d;\uN)$ determined by the
map $\Phi_\alpha$ on $\cO(d;\uN)$ (we note that the deformation
parameter is $\alpha-1$, not $\alpha$):
\begin{multline}
\label{comm}
{}[F,G]_{B,\alpha} :=
\Phi_\alpha^{-1}([\Phi_\alpha(F),\Phi_\alpha(G)]) =
\sum\limits_{1\leq i,j\leq d} B_{ij}\Phi_\alpha^{-1}(\del_i
\Phi_\alpha(F) \cdot \del_j
\Phi_\alpha(G)) \\
=\sum\limits_{1\leq i,j\leq d} B_{ij}\Phi_\alpha^{-1}(\del_i
\Phi_\alpha(F))\cdot \Phi_\alpha^{-1}(\del_j \Phi_\alpha(G)) =
 \sum\limits_{1\leq i,j\leq d} B_{ij} D_{\alpha,i}F\cdot D_{\alpha,j}G,
\end{multline}
because the map $\Phi_\alpha$ on $\cO(d;\uN)$ for $\alpha\neq 0$
preserves the (associative and commutative) multiplication of
functions.\footnote{Any automorphism of the \textit{space}
$\cO(d;\uN)$ produces a deformed bracket, but the second equality in
\eqref{comm} is due the fact that $\Phi_\alpha^{-1}(\del_i
\Phi_\alpha(F) \cdot \del_j \Phi_\alpha(G)) =
\Phi_\alpha^{-1}(\del_i \Phi_\alpha(F))\cdot \Phi_\alpha^{-1}(\del_j
\Phi_\alpha(G))$.}

We now consider the Lie algebra with the bracket \eqref{comm} for any
$\alpha$. Because we obtained this bracket from a trivial
deformation (for $\alpha\neq 0$), the obtained Lie algebra is
isomorphic to the initial Lie algebra $\fpo_B(d;\uN)$. To what is the
Lie algebra for $\alpha=0$ isomorphic?

Under the isomorphism between $\cO(d;\uN)$ and
$\cO(2d;(1,\dots,1,N_1-1,\dots,N_d-1))$ given by the formula
\begin{equation}\label{}
x_1^{(r_1)}\dots x_d^{(r_d)}\longleftrightarrow y_1^{(r_1\bmod
p)}\dots y_d^{(r_d\bmod p)}
y_{d+1}^{\left(\left[\frac{r_1}{p}\right]\right)}\dots
y_{2d}^{\left(\left[\frac{r_d}{p}\right]\right)},
\end{equation}
the operator $D_{0,i}$ on $\cO(d;\uN)$ becomes $\del_i$ on
$\cO(2d;(1,\dots,1,N_1-1,\dots,N_d-1))$. Hence, the Lie algebra given by
commutation relation (\ref{comm}) with $\alpha = 0$ is isomorphic to
\begin{equation}\label{Poisdefs}
\begin{aligned}
&\fpo_B(d;(1,\dots,1))\otimes \cO(d;(N_1-1,\dots,N_d-1))\\
&\qquad\simeq\fpo_B(d;(1,\dots,1))\otimes \cO(1;(N_1+\dots+N_d-d))\\
&\qquad\simeq\fpo_B(d;(1,\dots,1))\otimes\cO(1;(1))^{\otimes N_1+\dots+N_d-d}.
\end{aligned}
\end{equation}

We see from \eqref{Poisdefs} that all Poisson algebras with the same
number of indeterminates, the same $\sum N_i$, and bilinear forms
$B$ equivalent over the ground field are deforms of one Lie algebra.

Conjecturally, the statement
\begin{equation}\label{minip}
\begin{minipage}[l]{12cm}
``any vectorial Lie algebra $X(k;\uN)$ is a deform of the tensor
product \\ $X(k;\uN_s)\otimes\cO$, where $\cO:=\cO(u;\tilde\uN)$ with
an appropriate $\tilde\uN$"
\end{minipage}
\end{equation}
holds whenever the space of the Lie algebra $X(k;\uN)$ can be
identified with $\cO$, or the direct sum of several copies of $\cO$
each endowed with its extra structure of (associative)
multiplication, and the bracket can be defined using only
distinguished derivatives, associative and commutative
multiplication of functions, and linear operations, e.g., for
$\fvect$ and $\fk$ (cf.~the proof of Theorem 5 in\cite{Dzh}).

\section{The Jurman algebra is a semitrivial deform}\label{secJur}

\ssec{The Jurman algebra}
Jurman introduced a Lie algebra over $\Fee_2=\{0,1\}$ that seemed to
have no analog over fields $\Kee$ of characteristic $p\neq 2$
(see \cite{Ju}) until its interpretation in \cite{GJu}. Jurman
constructed this algebra by, in a sense, doubling the
\textit{Zassenhaus algebra}, which is the derived of the
\textit{Witt algebra} $\fvect(1;\uN)$. Jurman therefore called his
algebra the \textit{Bi-Zassenhaus algebra} denoted by $B(g,h)$. But
the letter $B$ is overused, and we wish to emphasize properties of
the Lie algebra $B(g,h)$ that differ from the properties in which
Jurman was interested. We therefore designate this algebra $\fj(g,h)$
in honor of Jurman. The following description (see \cite{Ju}) allows
extending the ground field and considering $\fj(g,h)$ over $\Kee$.

Let $g\geq 2$ and $h\geq 1$ be integers. Let $\eta= 2^g -1$ and
$\varkappa=2^{g+h}\geq 8$. Taking the elements
\begin{equation}\label{basisJur}
\bigl\{Y_{j}(t)\mid t\in\{0,1\},\;j\in \{-1,0,\dots,\varkappa-3\}\bigr\}
\end{equation}
as a basis in $\fj(g,h)$, Jurman defined the bracket by setting
\begin{equation}\label{brackJur}
[Y_{i}(s),Y_{j}(t)] = b_{s,t}^{i,j}\;Y_{i+j+st(1-\eta)}(s+t),
\end{equation}
where (in the next formula, the binomial coefficients and their sum
are considered modulo $2$ and meaningless expressions are taken to be
$0$; see Example \ref{EX} for further elucidations of the meaning of
the binomial coefficient for $s=t=1$)
{\small\begin{equation}\label{structcoefJur}
b_{s,t}^{i,j}=\begin{cases}
%{\!\textstyle
\binom{i + j + st (2-\eta)}{i + 1}+
%{\textstyle
\binom{i + j + st (2-\eta)}{j + 1}  %}\\
&%\phantom{0\kern30mm}
\text{if $-1\leq i + j +st (2-\eta)\leq \varkappa -3$},\\
0%\kern30mm
&\text{otherwise}.\end{cases} %\right.
\end{equation}
}

\sssbegin{Example}\label{EX} Let $(g, h)=(2, 1)$. We have
$b_{1,1}^{i,-1}=\binom{i-2}{i+1}+\binom{i-2}0$. The first summand is
meaningless for any $i$ and should be understood as a $0$. The
second summand makes sense for $i\ge 2$ when it is equal to 1.

We have $b_{1,1}^{i,0}=\binom{i-1}{i+1}+\binom{i-1}1$. The first
summand makes no sense for any $i$, and the second summand makes no
sense for $i=-1, 0, 1$. Each of these meaningless binomial
coefficients should be understood as $0$. If $i>1$, then
$\binom{i-1}1\equiv i-1\pmod 2$.

We have $b_{1,1}^{i,1}=\binom{i}{i+1}+\binom{i}2$ with the first
summand always meaningless (hence equal to $0$) and the second
summand equal to $0$ for $i<2$. \end{Example}

\ssec{The Jurman algebra $\fj(g,h)$ as a deform of $\fh_\Pi'(2; (g,
h+1))$} To interpret the Jurman algebra $\fj(g,h)$ somehow, we
compare it with a known simple Lie algebra. The most plausible
comparison candidates can be found in \cite{LeP}, where all possible
versions of Poisson Lie algebras were described in characteristic 2,
as well as  Lie (sub)algebras of Hamiltonian vector fields. We
realize the Poisson Lie algebra $\fpo_\Pi(2;\underline N)$ by
generating functions (divided powers) in the two indeterminates $p$
and $q$ with the bracket
\begin{equation}\label{bracketPo}
\{F,G\}=\frac{\del F}{\del p}\frac{\del G}{\del q}+
\frac{\del F}{\del q}\frac{\del G}{\del p}\quad
\text{for any }F,G\in\cO(2;\underline N),
\end{equation}
where $\del_p$ and $\del_q$ are distinguished partial derivatives.

We consider the Lie algebra of Hamiltonian vector fields $\fh_\Pi(2;
\underline N)=\fpo_\Pi(2; \underline N)/\Kee\cdot 1$ and its derived
$\fh_\Pi'(2; \underline N)$. We keep expressing the elements of $\fh_\Pi$ and
$\fh_\Pi'$ via generating functions having in mind, by abuse of
notation, their classes modulo the center of $\fpo_\Pi$.

We recall (see \cite{LSh1}) that the \textit{Weisfeiler filtrations}
were initially used to describe infinite-dimensional vectorial Lie
(super)algebras $\cL$ by selecting a maximal subalgebra $\cL_{0}$ of
finite codimension. Dealing with finite-dimensional algebras, we can
confine ourselves to maximal subalgebras of \textit{least} or
``almost least" codimension. Let $\cL_{-1}$ be the minimal
$\cL_{0}$-invariant subspace strictly containing $\cL_{0}$. For
$i\geq 1$, we set
\begin{equation}\label{1.3}
\cL_{-i-1}=[\cL_{-1}, \cL_{-i}]+\cL_{-i}\quad\text{and}\quad
\cL_i =\{D\in \cL_{i-1}\mid [D,\cL_{-1}]\subset\cL_{i-1}\}.
\end{equation}
We thus obtain a filtration,
\begin{equation}\label{1.2}
\cL= \cL_{-d}\supset \cL_{-d+1}\supset \dots \supset
\cL_{0}\supset \cL_{1}\supset \ldots.
\end{equation}
The $d$ in (\ref{1.2}) is called the \textit{depth} of $\cL$ and of
the associated graded Lie superalgebra
${\fg=\mathop{\oplus}\limits_{-d\leq i}\fg_i}$, where
$\fg_{i}=\cL_{i}/\cL_{i+1}$.

We let $\mathcal L$ denote $\fj(g,h)$ when considered with a
\textit{Weisfeiler filtration}. Eqs. \eqref{brackJur} and
\eqref{structcoefJur} imply that assuming $g,h=\infty$, we have
only one maximal subalgebra of finite codimension:
\[
\mathcal L_0=\Span(Y_i(0),Y_j(1)\mid i,j\ge 0).
\]
Its  maximality follows from table \eqref{basisfirst5}. The
Weisfeiler filtration corresponding to the pair
$(\mathcal L,\mathcal L_0)$ is
\begin{equation}\label{Wfil}
\cL= \cL_{-1}\supset \mathcal L_0\supset \mathcal L_1\supset
\mathcal L_2\dots,\quad\text{where }\mathcal L_{i+1}=\{X\in \mathcal
L_i\mid [\mathcal L,X]\subset \mathcal L_i\}.
\end{equation}
Let $\gr \fj(g,h)=\mathop{\oplus} \fg_i$, where
$\fg_i=\mathcal L_i/\mathcal L_{i+1}$ for $i\geq -1$.

\sssbegin{Proposition}\label{juriso} We have $\gr \fj(g,h)\cong
\fh_\Pi'(2;(g,h+1))$.
\end{Proposition}

\begin{proof} For brevity, we set $\fh:=\fh_\Pi'(2;(g,h+1))$. We first
note that every element of the Cartan prolong is uniquely determined
by its brackets with the elements of the $(-1)$st component. In
particular, any element $X=p^{(\beta)}q^{(\gamma)}\in\fh$ is uniquely
determined by the conditions
\begin{equation}\label{drh}
\begin{alignedat}{3}
&(\ad_p)^\gamma(\ad_q)^{\beta-1}(X)=p,&\quad
&(\ad_p)^{\gamma-1}(\ad_q)^{\beta}(X)=q&\quad
&\text{for }\beta\gamma> 0,\\
&(\ad_q)^{\beta-1}(X)=p,&\quad
&\ad_p(X)=0&\quad&\text{for }\beta>1,\;\gamma=0,\\
&(\ad_p)^{\gamma-1}(X)=q,&\quad
&\ad_q(X)=0&\quad&\text{for }\beta=0,\;\gamma>1.
\end{alignedat}
\end{equation}

We now pass to $\fg:=\gr \fj(g,h)$. Let $\bar X$ be the image of an
arbitrary element $X\in \fj(p,q)$ in $\fg$. The definition of a
filtration implies that $\dim \fg_{-1}=2$ and
$\fg_{-1}=\Span(\overline{Y_{-1}(0)}, \overline{Y_{-1}(1)})$. We
identify
\[
\overline{Y_{-1}(0)}\longleftrightarrow q,\quad
\overline{Y_{-1}(1)}\longleftrightarrow p.
\]

Let $0\leq \alpha\leq 2^h-1$ and $0\le\beta\le \eta=2^g-1$. Our next
goal is to establish the correspondence
\begin{equation}\label{iso1}
\overline{Y_i(s)}\longleftrightarrow p^{(\beta)}q^{(2\alpha+1-s)},\quad
\text{where }s=0,1\text{ and }i=\alpha(\eta+1)-1-s+\beta.
\end{equation}
For manual computations, it is more convenient to consider the two
cases $s=0,1$ separately:
\begin{equation}\label{iso2}
\begin{alignedat}{2}
&\overline{Y_a(1)}\longleftrightarrow p^{(\beta)}q^{(2\alpha)}&\quad
&\text{for }a=\alpha(\eta+1)-2+\beta,\\
&\overline{Y_b(0)}\longleftrightarrow p^{(\beta)}q^{(2\alpha+1)}&\quad
&\text{for }b=\alpha(\eta+1)-1+\beta.
\end{alignedat}
\end{equation}

Let $A_{\gamma,\delta}:=\left(\ad_{Y_{-1}(0)}\right)^{\gamma}
\left(\ad_{Y_{-1}(1)}\right)^{\delta}$. Clearly, the image $\bar X$
of $X\in \fj(p,q)$ belongs to $\fg_k$ if and only if there exist
$\gamma$ and $\delta$ such that $\gamma+\delta=k+1$ and $A_{\gamma,
\delta}(X)\notin \cal L_0$ while for $\gamma$ and $\delta$ such that
$\gamma+\delta<k+1$, we have $A_{\gamma, \delta}(X)\in \cal L_0$.

We now consider the brackets in the Lie algebra $\fj(p,q)$:
\begin{equation}\label{drj}
\begin{aligned}
&[Y_{-1}(0), Y_j(s)]=\begin{cases} Y_{j-1}(s)&\text{if }j\ge 0,\\
0&\text{if } j=-1, \end{cases} \\
&[Y_{-1}(1), Y_j(0)]=\begin{cases} Y_{j-1}(1)&\text{if }j\ge 0,\\
0&\text{if } j=-1, \end{cases} \\
&[Y_{-1}(1), Y_j(1)]=\begin{cases} 0&\text{if }j<\eta-1,\\
Y_{j-\eta}(0)&\text{if } j\ge \eta-1. \end{cases}
\end{aligned}
\end{equation}

Equations \eqref{drj} for the elements $X\in \fj(p,q)$ of the form
$X=Y_{\beta-2}(1)$ for $2\le \beta\le\eta$ imply that
\begin{equation}\label{which}
A_{\beta-1,0}(X)=Y_{-1}(1)\longleftrightarrow p, \quad A_{0,1}(X)=0.
\end{equation}
Expressions \eqref{which} mean that $X\in \fg_{\beta-2}$, and the
element $\bar X$ corresponds to $p^{(\beta)}\in\fh$. We have thus
obtained the first correspondence in \eqref{iso2} for $\alpha=0$.

Similarly, for $X=Y_{\beta-1}(0)$, where $1\le \beta\le\eta$, we
have
\[
\kern-4pt A_{\beta,0}(X)=Y_{-1}(0)\longleftrightarrow q \text{ and } A_{\beta-1,
1}(X)=Y_{-1}(1) \Longrightarrow X\in \fg_{\beta-1} \text{ and }
X\longleftrightarrow p^{\beta}q\in\fh,
\]
implying the second correspondence in eq.~\eqref{iso2} for $\alpha=0$.

Equations \eqref{drj} also imply that
\begin{equation}\label{y12}
\left(\ad_{Y_{-1}(1)}\right)^2\left(Y_j(s)\right)=\begin{cases}
Y_{j-\eta-1}(1)&\text{if }j\ge \eta,\\
0&\text{if }j<\eta.\end{cases}
\end{equation}
Therefore, for any $\alpha>0$ and $ \beta>0$ and
$X=Y_{\alpha(\eta+1)+\beta-2}(1)$, we have
\[
\renewcommand{\arraystretch}{1.2}
\begin{array}{l}
\left(\ad_{Y_{-1}(0)}\right)^{\beta-1}\left(\ad^2_{Y_{-1}(1)}\right)^\alpha(X)=
A_{\beta-1,2\alpha}(X)=Y_{-1}(1)\longleftrightarrow p,\\
\left(\ad_{Y_{-1}(0)}\right)^{\beta}\left(\ad^2_{Y_{-1}(1)}\right)^{\alpha-1}\ad_{Y_{-1}(1)}(X)=
A_{\beta,2\alpha-1}(X)=Y_{-1}(0)\longleftrightarrow q,
\end{array}
\]
implying the correspondence $\bar X\longleftrightarrow
p^{(\beta)}q^{(2\alpha)}$. This provides the first correspondence
in eq.~\eqref{iso2} for the case where $\alpha>0$ and $\beta>0$.
We obtain the second correspondence in eq.~\eqref{iso2} for the case where
$\alpha>0$ and $\beta>0$ absolutely analogously.

It remains to consider the case where $\alpha>0$ and $\beta=0$. Let
$X=Y_{\alpha(\eta+1)-1}(0)$. Then
\[
\text{Since $A_{0,2\alpha}(X)=Y_{-1}(0)\longleftrightarrow q$, it
follows that $\bar X\in\fg_{2\alpha-1}$}.
\]
But
\[
[Y_{-1}(0),X]=Y_{\alpha(\eta+1)-2}(0)=Y_{(\alpha-1)(\eta+1)+(\eta-1)}(0)
\longleftrightarrow p^{(\eta)}q^{(2\alpha-1)}\in\fg_{\eta+2\alpha-3}.
\]
Because $\eta\ge 3$, it follows that $2\alpha-1<\eta+2\alpha-3$, and hence
$[\overline{Y_{-1}(0)},\bar X]=0$, i.e., the element $X$ corresponds to
$q^{(2\alpha+1)}$. This provides the second correspondence in
eq.~\eqref{iso2} for $\alpha>0$ and $\beta=0$. The first correspondence
in eq.~\eqref{iso2} for this case is obtained similarly.

We see that the maximal power of $p$ is equal to $\eta=2^g-1$, and
hence $\underline N(p)=g$. Because $2\alpha+1\leq 2^{h+1}-1$, it follows
that $\underline N(q)=h+1$.\end{proof}

Accordingly, the basis elements of the components of the first five
degrees are as follows:
{\tiny
\begin{equation}\label{basisfirst5}
\renewcommand{\arraystretch}{1.4}\begin{array}{|c|c|c|c|c}\hline
\fg_{-1}&\fg_0&\fg_1& \fg_2 & \fg_3\\
\hline p\longleftrightarrow \overline{Y_{-1}(1)}&
p^{(2)}\longleftrightarrow \overline{Y_{0}(1)} &
p^{(3)}\longleftrightarrow \overline{Y_{1}(1)}
&p^{(4)}\longleftrightarrow \overline{Y_{2}(1)} &
p^{(5)}\longleftrightarrow \overline{Y_{3}(1)}\\
q\longleftrightarrow \overline{Y_{-1}(0)}& pq\longleftrightarrow
\overline{Y_{0}(0)}
& p^{(2)}q\longleftrightarrow
\overline{Y_{1}(0)}
& p^{(3)}q\longleftrightarrow
\overline{Y_{2}(0)}&
p^{(4)}q\longleftrightarrow \overline{Y_{3}(0)}\\
&q^{(2)}\longleftrightarrow \overline{Y_{\eta-1}(1)}&
pq^{(2)}\longleftrightarrow \overline{Y_{\eta}(1)} &
p^{(2)}q^{(2)}\longleftrightarrow \overline{Y_{\eta+1}(1)}
&p^{(3)}q^{(2)}\longleftrightarrow \overline{Y_{\eta+2}(1)}\\
&&q^{(3)}\longleftrightarrow \overline{Y_{\eta}(0)}
&pq^{(3)}\longleftrightarrow \overline{Y_{\eta +1}(0)}
&p^{(2)}q^{(3)}\longleftrightarrow \overline{Y_{\eta+2}(0)}\\
&&&q^{(4)}\longleftrightarrow \overline{Y_{2\eta}(1)}
&pq^{(4)}\longleftrightarrow Y_{2\eta+1}(1) \\
&&&&q^{(5)}\longleftrightarrow \overline{Y_{2\eta+1}(0)}\\
\hline\end{array}
\end{equation}}

We let $\{\,\cdot\,,\,\cdot\,\}$ and $[\,\cdot\,,\,\cdot\,]$ denote
the respective brackets in $\fg=\fh_\Pi'(2;(g,h+1))$ and $\fj(g,h)$.
Expressing the $\overline{Y_i(s)}$ in terms of monomials in $p$ and
$q$, we see that for the simplest case $g=h+1$, the \textit{Jurman
cocycle} $c$ (which deforms $\fh_\Pi'(2; (g,g))$ into the Jurman
algebra) is (for any $F,G\in\cO(2; (g,g))$), as direct calculations
show,
\begin{equation}\label{Jurbra}
[F,G]=\{F,G\}+m_c(F,G),\quad\text{where }
c=\mathop{\sum}\limits_{m<n}p^{(\eta)}q^{(m+n-3)}\otimes
d(q^{(m)})\wedge d(q^{(n)}),
\end{equation}
and $m_c(F,G)$ (see eq.~\eqref{JurMap}) is the map corresponding to
the cocycle $c$.

We call all other cocycles (which do not deform $\fh_\Pi'(2;(g,g))$
to the Jurman algebra) \textit{non-Jurman} cocycles.

With respect to the pair of operators $(\deg_p(\cdot)-1,\deg_q(\cdot)-1)$
the weight of the cocycle $F\otimes d(G_1) \wedge \dots \wedge d(G_n)$
for any $F,G_1,\dots,G_n\in \fh_\Pi'(2; (g,g))$
is equal to
\begin{equation}\label{Zgra}
(\deg_p(F)-1-\sum (\deg_p(G_i)-1),\quad\deg_q(F)-1-\sum(\deg_q(G_i)-1)).
\end{equation}
Hence, the Jurman cocycle has the weight $(2^g,-2)$. By the symmetry
$p\longleftrightarrow q$, there is \textit{another Jurman cocycle} of
weight $(-2, 2^g)$ leading to an isomorphic Jurman algebra.

If $g\neq h+1$, then there is no symmetry $p\longleftrightarrow q$, but
there still is \textit{another Jurman cocycle} making $\fh_\Pi'(2;
(g, h+1))$ into $\fj(h+1, g-1)$. It is of the form
\begin{equation}\label{Jurbra1}
c=\mathop{\sum}\limits_{m<n}q^{(\theta)}p^{(m+n-3)}\otimes
d(p^{(m)})\wedge d(p^{(n)}),\text{~~where $\theta=2^{h+1}-1$.}
\end{equation}

\paragraph{Remark}In characteristic $p>2$,  most of the cocycles
representing classes of $H^2(\fg;\fg)$ are not integrable for the
simple vectorial Lie algebras $\fg$ (see \cite{Dzh}). If $p=2$, we
do not know any simple Lie algebra $\fg$ with a nonintegrable
cocycle representing a class of $H^2(\fg;\fg)$.

\sssbegin{Lemma}[For $(g,h)=(2,1), (2,2)$; it is a conjecture for
generic values of $(g,h)$]\label{conjlemma} Any linear combination
of cocycles representing classes of $H^2(\fg;\fg)$ for
$\fg=\fh_\Pi'(2;(g,h+1))$  can be integrated to a global deform
(Cf.~$\fg=\fh_I(2;(g,h+1))$; see Lemma $\ref{nonlin}$).
\end{Lemma}

%\begin{proof} The lemma is proved by computer-aided study for
%$(g,h)=(2,1), (2,2)$.\end{proof}

For $g+h=g'+h'=K$, the Jurman algebras $\fj(g,h)$ and $\fj(g',h')$
regarded as $\Zee/2$-graded Lie algebras $\fj=\fj_\ev\oplus\fj_\od$
with $\fj_\ev$ spanned by the $Y_i(0)$ for all $i$ have these even
parts isomorphic, and the odd parts, as modules over the even part,
are also isomorphic. This is clear from eqs.~\eqref{brackJur} and
\eqref{structcoefJur}. We note that the brackets of two odd elements
given by Jurman's cocycles can be united into one bracket depending
on as many parameters as there are partitions $K=g+h$ with $g\geq 2$
and $h\geq 1$. To see this, we consider the brackets of two ``odd"
elements and one ``even" element and also consider the brackets of
three ``odd" elements; the statement is obvious in both cases. The
obtained bracket depends linearly on all $K$ parameters.

\ssbegin{Proposition}\label{JurmanIso} The Jurman algebra $\fj(g,h)$
is isomorphic to $\fh:=\fh_\Pi'(2;(g,h+1))$. For $0\leq k <2^{h+1}$
and $0\leq l<2^{g-1}$, the isomorphism is given by the maps
\begin{equation}\label{AL0}
Y_{2^{g-1}k+l-1}(0) \longleftrightarrow Y_{2^{g-1}k+l-1} =
p^{(2^{g-1}+l)}q^{(k)} + (k+1)p^{(l)}q^{(k+1)},
\end{equation}
\begin{equation}\label{AL1} Y_{2^{g-1}k+l-2}(1) \longleftrightarrow
Z_{2^{g-1}k+l-2} =\begin{cases} p^{(2^{g-1}+l)}q^{(k-1)} +
(k+1)p^{(l)}q^{(k)}&\text{if }k>0,\\
p^{(l)}&\text{if }k=0.\end{cases} \end{equation}
Assuming that $q^{(m)}=0$ for $m<0$, we can express $Z$ in \eqref{AL1}
uniformly:
\begin{equation}\label{AL2}
Z_{2^{g-1}k+l-2}=p^{(2^{g-1}+l)}q^{(k-1)}+(k+1)p^{(l)}q^{(k)}\quad
\text{for all values }0\leq k < 2^h.
\end{equation}
\end{Proposition}

\begin{proof} Taking into account that $(k+1)\mod 2=0$ or 1, we can
unite correspondences \eqref{AL0} and \eqref{AL1} by considering the
1-parameter family of maps (we note that $\ad_{p^{(\mu+1)}}$ is a
derivation of $\fh$ such that $\ad_{p^{(\mu+1)}}^2=0$):
\[
\Phi_t:\fh\longrightarrow \fh,\quad F\mapsto
F+tp^{(\mu)}\frac{\del F}{\del q}=F+t\ad_{p^{(\mu+1)}}(F),\quad
\text{where }\mu=2^{g-1}.
\]
Let $m_c(F,G)$ be the map \eqref{Jmap} corresponding to the cocycle $c$
(see \eqref{Jurbra}). Direct computations show that $\Phi_{\sqrt{t}}$
is an isomorphism of Lie algebras because
\[
\{\Phi_t(F), \Phi_t(G)\}=\Phi_t(\{F,G\})+t^2m_c(F,G). \qed
\]
\noqed\end{proof}

\sssec{Deforms of $\fh_\Pi'(2;(g,h+1))$ for the smallest values of
$(g,h)$} If $g=h+1$, it clearly suffices to consider only cocycles
of nonnegative weight because of the symmetry $p\longleftrightarrow
q$.

Let $(g, h)=(2,1)$. Here are cocycles representing a basis of the space
$H^2(\fg;\fg)$. The Jurman cocycle $c:=c_{2^g,-2}$ in~\eqref{Jurbra} is
$c_{4,-2}$ from our list \eqref{gh=21cocycles} below.

\ssbegin{Proposition}\label{Jcocycles} Here, $F,G\in\fh_\Pi'(2;(g,h+1))$
are arbitrary, $\{\,\cdot\,,\,\cdot\,\}$ is the Poisson bracket of
functions generating $\fh_\Pi(2m;\uN)$, and $\hbar\in\Kee$.

\textbf{The following cocycles are semitrivial}:

1. The Jurman cocycle $c_{2^{g},-2}$ (see \eqref{Jurbra}) represents
the map
\begin{equation}\label{Jmap}
p^{(\eta)}(\del_q F\cdot\del_q^2 G+\del_q^2 F\cdot\del_q G),\quad
\text{where }\eta=2^{g}-1.
\end{equation}
For $(g,h)=(2,1)$, the cocycle $c_{4,-2}$ (see \eqref{gh=21cocycles})
represents the map
\begin{equation}\label{JurMap}
m_{4,-2}(F,G)=p^{(3)}(\del_q F\cdot\del_q^2 G+\del_q^2 F\cdot\del_q G).
\end{equation}

2. The cocycle $c_{0,-4}$ (see \eqref{gh=21cocycles}) represents the
map whose shape is independent of $(g,h)$:
\begin{equation}
m_{0,-4}(F,G)= \del_p\del_q^2 F \cdot \del_q^3 G + \del_q^3 F\cdot
\del_p\del_q^2 G = \{\del_q^2 F, \del_q^2 G\}.
\end{equation}

3. The cocycle $c_{0,-2}$ (see \eqref{gh=21cocycles}) is equivalent
to the cochain that represents the map whose shape is independent
of $(g,h)$:
\begin{equation}
m_{0,-2}(F,G)= \{\del_q F, \del_q G\}.
\end{equation}

\textbf{The following cocycles are nontrivial}:

4. The cocycle $c_{2,0}$ (see \eqref{gh=21cocycles}) is equivalent
to the cochain representing the map
\begin{equation}
m_{2,0}(F,G)= p^{(\eta)}(\del_q F\cdot \del_p^2 G + \del_p^2 F\cdot
\del_q G),\quad\text{where }\eta=2^{g}-1,
\end{equation}
which yields one of the filtered deforms (see \cite{Dzh}).

5. The cocycle $c_{-2,-2}$ (see \eqref{gh=21cocycles}) is inherited
from the quantization of the Poisson Lie algebra $\fpo_\Pi(2;(a,a))$
being the linear in the Planck constant part of the cocycle
restricted to the subquotient $\fh_\Pi'$ of $\fpo_\Pi$. The
deformation turns $\fh_\Pi'(2;(a,a))$ into $\fpsl(2^a)$ for any $a$.
\end{Proposition}

The index of each cocycle is equal to its weight (further on, to
save trees, we give full expression of only the cocycles with short
expressions; the lexicographic order of summands adding up to the
cocycle makes it possible to distinguish cocycles by looking at the
pieces displayed; if these pieces are insufficient to interpret them,
then see the \TeX file in arXiv and uncomment the hidden terms):
{\tiny
\begin{equation}\label{gh=21cocycles}
\begin{aligned}
&\begin{aligned}
c_{4,-2}={}&p{}^{(3)}\otimes d(q)\wedge d(q{}^{(2)}) + p{}^{(3)}\, q
\otimes d(q)\wedge d(q{}^{(3)}) + p{}^{(3)}\, q{}^{(2)}\otimes
d(q{}^{(2)})\wedge d(q{}^{(3)}),
\end{aligned}\\[1mm]
&\begin{aligned}
c_{0,-4}={}&p\otimes d(p\, q{}^{(2)})\wedge d(p\, q{}^{(3)}) +
p\otimes d(q{}^{(3)})\wedge d(p{}^{(2)}\, q{}^{(2)}) + q\otimes
d(q{}^{(3)})\wedge d(p\, q{}^{(3)})+\dots\\
%&{}+p{}^{(2)}\otimes d(p\, q{}^{(2)})\wedge d(p{}^{(2)}\, q{}^{(3)}) +
%p{}^{(2)}\otimes d(q{}^{(3)})\wedge d(p{}^{(3)}\, q{}^{(2)})\\
%&{}+p\,q\otimes d(q{}^{(3)})\wedge d(p{}^{(2)}\,
%q{}^{(3)})+p{}^{(3)}\otimes d(p{}^{(2)}\, q{}^{(2)})\wedge d(p{}^{(2)}\,
%q{}^{(3)})\\
%&{}+p{}^{(3)}\otimes d(p\, q{}^{(3)})\wedge d(p{}^{(3)}\,
%q{}^{(2)})+p{}^{(2)}\, q\otimes d(p\, q{}^{(3)})\wedge d(p{}^{(2)}\,
%q{}^{(3)}),
\end{aligned}\\[1mm]
&\begin{aligned}
c_{2,0}={}&p^{(2)}\otimes d(p)\wedge d(q) + p\, q{}^{(2)}\otimes
d(q)\wedge d(q{}^{(2)}) +
p{}^{(3)}\, q\otimes d(q)\wedge d(p{}^{(2)}\, q)\\
&{}+p{}^{(3)}\, q{}^{(2)}\otimes d(p)\wedge d(p\, q{}^{(3)}) +
p{}^{(3)}\, q{}^{(2)}\otimes d(q{}^{(2)})\wedge d(p{}^{(2)}\, q) +
p{}^{(2)}\, q{}^{(3)}\otimes d(q)\wedge d(p\, q{}^{(3)}),
\end{aligned}\\[1mm]
&\begin{aligned}
c_{0, -2}={}&p\otimes d(p)\wedge d(p\, q{}^{(3)})+ p\otimes d(p\,
q)\wedge d(p\, q{}^{(2)}) +p\otimes d(q{}^{(2)})\wedge d(p{}^{(2)}\,q)\\
&{}+q\otimes d(q)\wedge d(p\, q{}^{(3)}) +
q\otimes d(p\, q)\wedge d(q{}^{(3)}) +q\otimes d(q{}^{(2)})\wedge
d(p\, q{}^{(2)})+\dots\\
%&{}+p{}^{(2)}\otimes d(p)\wedge d(p{}^{(2)}\,
%q{}^{(3)})+p{}^{(2)}\otimes d(p\, q)\wedge d(p{}^{(2)}\, q{}^{(2)})
%+p{}^{(2)}\otimes d(q{}^{(2)})\wedge d(p{}^{(3)}\, q)\\
%&{}+p\, q\otimes
%d(p\, q)\wedge d(p\, q{}^{(3)})+ q{}^{(2)}\otimes d(q{}^{(2)})\wedge d(p\, q{}^{(3)})+
%p{}^{(3)}\otimes d(p{}^{(2)})\wedge d(p{}^{(2)}\, q{}^{(3)})\\
%&{}+p{}^{(3)}\otimes d(p\, q)\wedge d(p{}^{(3)}\, q{}^{(2)})
%+p{}^{(3)}\otimes d(p{}^{(3)})\wedge d(p\,q{}^{(3)})
%+p{}^{(3)}\otimes d(p{}^{(2)}\, q)\wedge d(p{}^{(2)}\, q{}^{(2)})\\
%&{}+p{}^{(3)}\otimes d(p\, q{}^{(2)}) \wedge d(p{}^{(3)}\, q)
%+p{}^{(2)}\, q\otimes d(p\, q)\wedge d(p{}^{(2)}\, q{}^{(3)})\\
%&{}+p\,q{}^{(2)}\otimes d(q{}^{(2)})\wedge d(p{}^{(2)}\,
%q{}^{(3)})+p{}^{(3)}\, q\otimes d(p{}^{(2)}\,q)\wedge d(p{}^{(2)}\,
%q{}^{(3)})\\
%&{}+p{}^{(3)}\, q\otimes d(p{}^{(3)}\, q)\wedge d(p\,
%q{}^{(3)})+p{}^{(2)}\, q{}^{(2)}\otimes d(p\, q{}^{(2)})\wedge
%d(p{}^{(2)}\, q{}^{(3)})\\
%&{}+p{}^{(2)}\, q{}^{(2)}\otimes d(p{}^{(2)}\,
%q{}^{(2)})\wedge d(p\, q{}^{(3)}) +
%p{}^{(2)}\, q{}^{(3)}\otimes d(p\, q{}^{(3)})\wedge d(p{}^{(2)}\,
%q{}^{(3)})\\
%&{}+p{}^{(2)}\, q{}^{(3)}\otimes
%d(p\, q{}^{(3)})\wedge d(p{}^{(2)}\, q{}^{(3)}),
\end{aligned}\\[1mm]
&\begin{aligned}
c_{-2,-2}={}&p\otimes d(p\, q{}^{(2)})\wedge d(p{}^{(3)}\, q)+
q\otimes d(p\, q{}^{(2)})\wedge d(p{}^{(2)}\, q{}^{(2)})\\
&{}+q\otimes d(q{}^{(3)})\wedge d(p{}^{(3)}\, q )+
p{}^{(2)}\otimes d(p{}^{(3)}\, q)\wedge d(p{}^{(2)}\, q{}^{(2)}) +
p\, q\otimes d(p\, q{}^{(2)})\wedge d(p{}^{(3)}\, q{}^{(2)})+\dots\\
%&{}+p\, q
%\otimes d(p{}^{(3)}\, q)\wedge d(p\, q{}^{(3)}) + q{}^{(2)}\otimes
%d(p\, q{}^{(2)})\wedge d(p{}^{(2)}\, q{}^{(3)})\\
%&{}+q{}^{(2)}\otimes d(q{}^{(3)}) \wedge d(p{}^{(3)}\, q{}^{(2)}) +
%q{}^{(2)}\otimes d(p{}^{(2)}\, q{}^{(2)})\wedge d(p\,q{}^{(3)})\\
%&{}+p{}^{(3)}\otimes d(p{}^{(3)}\, q)\wedge d(p{}^{(3)}\, q{}^{(2)}) +
%p{}^{(2)}\, q\otimes d(p{}^{(3)}\, q)\wedge d(p{}^{(2)}\,q{}^{(3)})\\
%&{}+p{}^{(2)}\, q \otimes d(p{}^{(2)}\, q{}^{(2)})\wedge
%d(p{}^{(3)}\, q{}^{(2)}).
\end{aligned}
\end{aligned}
\end{equation}}

\begin{proof} The fact that the maps $m_w$, where $w$ is a weight, do
correspond to the cocycles $c_w$ as claimed is subject to a direct
verification. The idea is as follows. We see that the image of
$c_{4,-2}$ is always divisible by $p^{(3)}$ and that this image is
nonzero only if both arguments are polynomials (divided powers) in
$q$. Taking into account the weight of the cocycle and its
(anti)symmetry, we seek $m_c$ in the form
\[
p^{(3)}(A(F\del_q^3 G +
G\del_q^3 F) + B(\del_q F\cdot \del_q^2 G + \del_q^2 F \cdot\del_q
G)),\quad\text{where }A,B\in\Kee.
\]
It turns out that for $A=0$ and $B=1$, we obtain the desired. For
other cocycles, we seek the operators in the form $D_1F \cdot D_2G +
D_2F\cdot D_1G$ or $p^{(3)}(D_1F\cdot D_2G + D_2F\cdot D_1G)$, where
the $D_i$ are compositions of some derivations. We must check if any
of these operators $m_c$ in fact matches $c$. It could be that to
have a nice expression for $m_c$, we must replace the cocycle $c$
with a $\tilde c$ of the same cohomology class.

1. The semitriviality of the Jurman cocycle is explicitly proven for
arbitrary $(g,h)$ in Proposition \ref{JurmanIso}.

2. We consider the maps $\Phi_\hbar(F) = F + \hbar DF$, where $D =
\del_q^2$. Because $D^2=0$, it follows that the corresponding
deformed bracket produced by $c_{0,-4}$ is
\begin{equation}\label{sqrt}
{}\{F,G\}^{c_{0,-4}}_\hbar:=\{F,G\}^\Phi_{\sqrt{\hbar}}.
\end{equation}

3. In this case, although $D^2\neq 0$ for $D=\del_q$, the derivation $D$
is still nilpotent, and Theorem \ref{semiTr} is hence applicable here.

4. In this case, the deformed bracket is equivalent to
\begin{equation}
{}\{F,G\}_\hbar=(\del_p+\hbar p^{(3)}\del_q^2)F \cdot \del_q G+
\del_q F \cdot (\del_p + \hbar p^{(3)}\del_q^2)G.
\end{equation}

5. Let the Poisson Lie algebra be realized by the Poisson bracket on the
space of functions in $\vec p=(p_1,\dots,p_m)$ and $\vec q=(q_1,\dots,q_m)$.
We consider the deformation (over the ground field $\Fee=\Cee$ or $\Ree$,
physicists call it \textit{quantization}) that turns the Poisson Lie
algebra into the Lie algebra of differential operators with polynomial
coefficients (see sect.~1.4.7 in \cite{Fu}). The cocycle that determines
\textit{quantization} corresponds to the map
\begin{equation}\label{Qlin}
\cQ(F,G)=\mathop{\sum}\limits_{1\leq i\leq m}\frac{\del^2
F}{\del p_i^2}\frac{\del^2 G}{\del q_i^2}-\frac{\del^2 F}
{\del q_i^2}\frac{\del^2 G}{\del p_i^2}\quad\text{for any }F,G\in\Fee[p,q].
\end{equation}
Here, we encounter an analog of quantization over $\Fee=\Kee$ for
$\Char\Kee=2$. Let the coordinates of the shearing vector
corresponding to $\vec p$ be the same as those corresponding to
$\vec q$. Let $\hat{}: F\mapsto \hat F$ be the map that to any
monomial $F\in \cO(\vec p,\vec q;(\uN, \uN))$ ordered such that each
$p_i$ is to the left of all the $q_j$ for all $i$ and $j$ assigns a
differential operator obtained by the replacement $q_i\mapsto
\hbar\del_{p_i}$, where $\hbar\in\Kee$, for each $i$. All linear
operators in the finite-dimensional space $\cO(\vec p;\uN)$ are
differential, and so the deformed Lie algebra is isomorphic to
$\fgl(\cO(\vec p;\uN)) \simeq\fgl(2^{|\uN|})$, where $|\uN|=\sum
\uN_i$. Clearly, the same cocycle induces a deformation of
$\fh'(2m;\uN)$ into $\fpsl(2^{|\uN|})$. For any $\hbar\neq0$, the
deforms are obviously isomorphic (use rescaling, i.e., divide by
$\hbar$), and the commutator of differential operators is related to
the Poisson bracket as
\begin{equation}\label{hbar}
[\hat F,\hat G]=\{F,G\}_{P.b.} + O(\hbar)\quad\text{for any }
F,G\in \cO(\vec p,\vec q;(\uN, \uN)).
\end{equation}
For $m=1$, the weight of the cocycle part linear in $\hbar$ in the
right-hand side of eq.~\eqref{hbar} (up to a sign corresponding to
the interchange $p\leftrightarrow q$)  is precisely $-(2,2)$.
\end{proof}

For $\fg=\fh_\Pi'(2;(g,h+1))$, where $(g,h)=(2,1)$, the cocycles
considered in Proposition \ref{Jcocycles} represent a basis of
$H^2(\fg;\fg)$. Because $\dim H^2(\fg;\fg)$ increases together with
the coordinates of the shearing vector $(g,h)$, there are more
deformations to be interpreted in the general case. It seems
reasonable to switch attention from cocycles $c$ to maps $m_c$.
\textbf{Conjecturally}, all non-Jurman cocycles correspond to the
filtered deforms classified by Skryabin (see \cite{Sk}) or to the
quantization. This is so for $(g,h)=(2,1)$.

\section{What Kaplansky algebras are isomorphic to. Nonlinear
superizations}\label{secK} In 1981, Kaplansky described four types (in
fact, five: the dimensions of the two cases of the type-4 algebras
differ significantly) of simple Lie algebras for $p=2$ (see \cite{Kap2}).
He described them only by means of the multiplication table. We
interpret them in terms of familiar Lie algebras of vector fields.

Kaplansky defined the algebras in terms of $J$-\textit{systems}
resembling the notion of a root system. Over $\Fee_2$, a $J$-system
$\Gamma$ in the space $V$ with a symmetric inner product $B$ is a set
of nonzero vectors with the property that if $u, v\in \Gamma$ are
distinct and satisfy $B(u,v)=1$, then $u+v\in \Gamma$. Given any
$J$-system $\Gamma$, we construct a Lie algebra $\fg_\Gamma$ over
$\Fee_2$ with basis elements $e_u$ for every $u\in \Gamma$, and the
multiplication given by the expressions
\begin{equation}\label{Kapbracket}
{}[e_u, e_v]= \begin{cases}B(u,v)e_{u+v}&\text{for }u,v
\text{ distinct and }u+v\in \Gamma,\\
0&\text{for }u+v\not\in\Gamma\text{ or }u=v.\end{cases}
\end{equation}
We note that the second half of the lower property
in \eqref{Kapbracket} is automatically satisfied if the form $B$
is alternating. Each Kaplansky algebra $\text{Kap}_i(n)$, where
$i=1,2,3,4$, has the form $\fg_\Gamma$ for some $\Gamma$.

Any algebra defined over $\Fee_2$ can obviously be defined over $\Kee$
by extension of the ground field. In what follows, speaking of
Kaplansky algebras, we assume that such an extension is performed
unless otherwise specified.

\underline{$\text{Kap}_1(n)$}: For $n\geq 4$, let $\dim V=n$, and let
$V$ carry a nondegenerate and nonalternate inner product $B$. Let
$e^1,\dots,e^n$ be an orthonormal basis of $V$. For $\Gamma$, we take
all vectors in $V$ except $0$ and $e=e^1+\dots + e^n$, which can be
invariantly described as the unique element satisfying $B(e,y)=B(y,y)$
for all $y$.

Clearly, $B\sim I$. We recall that the brackets in
$\fl\fh_I(V; \uN):=(V,\fo_I'(V))_{*, \uN}$ ($\fl$ is for ``little")
and $\fh_I(V; \uN):=(V, \fo_I(V))_{*, \uN}$ are the same  (see \cite{LeP}):
\[
\{F,G\}_I:=\mathop{\sum}\limits_{1\leq i\leq n}\del_{z_i}F\cdot
\del_{z_i}G,\quad\text{where }F,G\in \Kee[z;\uN],
\]
only the stocks of generating functions of these Lie algebras differ.
We make the assignment
$e_u\longleftrightarrow\mathop{\prod}\limits_{1\leq i\leq n}(1+z_i)^{u_i}$.
(We note that Kaplansky considered monomials in $X_i:=1+z_i$ instead of
monomials in $z_i$.) We have $\text{Kap}_1(n)\simeq\fl\fh_I'(n)$ because
$\text{Kap}_1(n)$ does not contain $e$. In particular, we have an
interpretation of $\text{Kap}_1(4)$ sought, but not found,
in \cite{Ju, GJu}. Eick proved the isomorphism for $n=4$ (in different
terms) in \cite{Ei}.

\underline{$\text{Kap}_2(2m)$}: Let $\dim V=2m$, and let $V$ carry a
nondegenerate and alternate inner product $\Pi$. We take all nonzero
vectors in $V$. Kaplansky mentioned this algebra because it fits
into the approach he suggested although this algebra has analogs for
any characteristic\footnote{Kaplansky did not describe such algebras
explicitly. Here is a description for any $p>0$: Consider the
polynomial algebra in $y_i:=\exp(x_i)$ and set $\del_{x_i}y_j=
\delta_{ij} y_j$ and $(y_i)^p = \exp(px_i) = 1$. In the space
$\Kee[y]$, we introduce the Poisson bracket. Then $\text{Kap}_2(2m)$
is isomorphic to the quotient of the Poisson algebra modulo the
ideal of constants.} $p>0$, and we could hence have ignored it; it
is a filtered deform of $\fh_\Pi(2m;\uN_s)$. If we had ignored it,
then we would not have discovered a nonlinear $\Zee/2$-grading.

\underline{$\text{Kap}_3(n)=\fo_I'(n)$}, as Kaplansky observed (in
different terms). Kaplansky wrote ``the gaps (in the set of values
of $n=5$, 7, and $\geq 9$) avoid duplication."

\underline{$\text{Kap}_{4,a}(2m)$, where $a=0$ or 1}, is a temporary
notation, for lack of a better idea, for two similarly described and
equally mysterious algebras of quite different dimensions. In their
description, we need Arf invariants of quadratic forms. For a most
lucid definition of an Arf invariant, see \cite{Dye}. In
eq.~\eqref{defKap4A}, $a$ is the value of the Arf invariant (here,
$0$ or 1), and $B$ is short for ``Big" and is reminiscent of the
form $B$ (see eq.~\eqref{bilf}).

Let $\dim V=2m$, where $m\geq 3$, and let $Q$ be a nondegenerate
quadratic form on $V$. We set
\begin{equation}\label{defKap4A}
\begin{alignedat}{2}
&\text{Kap}_{4,a}(2m):=\fg_{\Gamma_a}(2m)&\quad
&\text{for }\Gamma_a=\{u\in V \mid Q(u)=1\},\text{ where Arf}(Q)=a,\\
&\text{Kap}_{4,B}(2m):=\fg_{\Gamma_B}(2m)&\quad
&\text{for }\Gamma_B=\{u\in V\},\end{alignedat}
\end{equation}
where the alternating bilinear form $B$ is given by the formula
\begin{equation}\label{bilf}
B(u,v)=Q(u+v)+Q(u)+Q(v).
\end{equation}
We note that several quadratic forms $Q$, nonequivalent and with
different values of the Arf invariant, in eq.~\eqref{bilf} can
produce the same bilinear form $B$. Observe that
\[
\text{Kap}_{4,a}(2m)\subset\text{Kap}_{2}(2m)\subset\text{Kap}_{4,B}(2m).
\]

\ssbegin{Proposition} 1. The Lie algebra \textit{Kap}${}_{4,B}(2m)$
is isomorphic to the algebra whose space is $\cO(2m;\uN_s)$ with
indeterminates $p_i$ and $q_i$ for $1\leq i\leq m$ and the bracket
\begin{equation}\label{brKap}
[f,g] = \mathop{\sum}\limits_{1\leq i\leq m} (1+ p_i)(1+
q_i)(\del_{p_i}f\cdot\del_{q_i}g + \del_{q_i}f\cdot\del_{p_i}g).
\end{equation}

2. The Lie algebra \textit{Kap}${}_{4,B}(2m)$ is isomorphic to a
deform of the Poisson algebra $\fpo_\Pi(2m;\uN_s)$ with the deformed
bracket
\begin{equation}\label{brKap-deformed}
[f,g]_\hbar = \mathop{\sum}\limits_{1\leq i\leq m} (1+\hbar
p_iq_i)(\del_{p_i}f\cdot\del_{q_i}g + \del_{q_i}f\cdot\del_{p_i}g)\quad
\text{for any }\hbar\neq 0
\end{equation}
and
\begin{equation}\label{Kap2}
\text{\textit{Kap}}{}_{4,B}(2m)\simeq\text{\textit{Kap}}{}_{2}(2m)\oplus\fc,
\end{equation}
where the center $\fc$ is generated by constant functions.
\end{Proposition}

\begin{proof} 1. The isomorphism is given as follows. We choose a
symplectic basis for the inner product $B$ in $V$. If $(u_1,\dots,u_{2m})$
are coordinates of a vector $u\in V$ in this basis, then
\[
e_u\longleftrightarrow f_{u}=(1+p_1)^{u_1}\dots
(1+p_m)^{u_m}(1+q_1)^{u_{m+1}}\dots (1+q_m)^{u_{2m}}.
\]

2. Clearly, \eqref{brKap} is a particular case of the bracket
\begin{equation}\label{brKap-deformed2}
[f,g]_\hbar = \mathop{\sum}\limits_{1\leq i\leq m} (1+\hbar'
p_i)(1+\hbar' q_i)(\del_{p_i}f\cdot\del_{q_i}g +
\del_{q_i}f\cdot\del_{p_i}g)\quad\text{with }\hbar'=1.
\end{equation}
Here, the part linear in $\hbar'$ describes a trivial deformation of
$\fpo_\Pi(2m;\uN_s)$ (as can be verified), and the quadratic part
corresponds to \eqref{brKap-deformed} with $\hbar=(\hbar')^2$. This
cocycle is nontrivial, as a computer-aided study shows.

The center is a direct summand because all weight spaces in
$\text{Kap}_{4,B}(2m)$ are 1-dimensional, and the weight of the space
generated by constants is $0$, but there are no two distinct weight
vectors of the same weight.
\end{proof}

\sssec{Kaplansky algebras $\text{Kap}_{4,B}(2m)$ and
$\text{Kap}_{4,a}(2m)$ in convenient indeterminates}
Examples of forms $Q_a$ with an Arf invariant equal to $a$ are
\begin{equation}\label{Q}
\begin{aligned}
&Q_{0}(u)=\mathop{\sum}\limits_{1\leq i\leq m} u_iu_{m+i},\\
&Q_{1}(u)=u_1^2+u_{m+1}^2+
\mathop{\sum}\limits_{1\leq i\leq m} u_iu_{m+i},
\end{aligned}
\end{equation}
We introduce operators $L_i$, where $i=1,\dots,2m$:
\[
L_i = \begin{cases} (1+p_i)\del_{p_i}&\text{if }1\leq i \leq m,\\
(1+q_{i-m})\del_{q_{i-m}}&\text{if }m+1\leq i\leq2m.\end{cases}
\]
Then $L_i f_u =u_if_u$. We set $\Delta=\mathop{\sum}
\limits_{1\leq i\leq m}L_iL_{i+m}$.

The subalgebras $\text{Kap}_{4,a}(2m)\subset\text{Kap}_{4,B}(2m)$
with bracket \eqref{brKap} are spanned by the nonzero elements
$f_u$ such that $Q_a(u)=1$. From the definition \eqref{defKap4A}, we
derive the conditions that single out the subalgebras
$\text{Kap}_{4,a}(2m)$ in $\text{Kap}_{4,B}(2m)$:
\begin{equation}\label{Qi}
\begin{alignedat}{2}
&f+\Delta f=0&\quad&\text{for }\text{Kap}_{4,0}(2m),\\
&f+(1+p_1)\del_{p_1}f+(1+q_1)\del_{q_1}f+\Delta f=0&\quad
&\text{for }\text{Kap}_{4,1}(2m).
\end{alignedat}
\end{equation}
The condition $\Delta f+f=0$ in \eqref{Qi} singles out the
eigenvectors of $\Delta$ with the eigenvalue $1$. But
\[
\Delta f_u =
\mathop{\sum}\limits_{1\leq i\leq m} u_iu_{i+m} f_u = Q_0(u)f_u,
\]
and this eigenspace is therefore spanned by all $f_u$ such that
$Q_0(u)=1$, which is exactly the image of $\text{Kap}_{4,0}$.

The case of $\text{Kap}_{4,1}$ is similar. For simplicity, we
respectively replace $L_1^2 f$ and $L_{m+1}^2 f$ with $L_1 f$ and
$L_{m+1} f$. This is possible because $L_i^2 = L_i$ and the $u_i$
only take values $0$ and $1$. Indeed,
\[
L_i^2 f_u = u_i^2 f_u = u_i f_u = L_i f_u.
\]

Kaplansky claimed (and we see that the claim obviously follows
from \eqref{Q}) that
\begin{equation}\label{dimKap4A}
\dim\fg_{\Gamma_a}=2^{m-1}(2^m- (-1)^a)=\begin{cases}2^{m-1}(2^m-1)
&\text{if Arf}(Q)=0,\\
2^{m-1}(2^m+1)&\text{if Arf}(Q)=1.\end{cases}
\end{equation}

We now study the structure of these algebras. It is more convenient
to pass to the coordinates $x_i:=(1+p_i)$ and $y_i:=(1+q_i)$. The
bracket \eqref{brKap} and operators \eqref{Qi} become
\begin{equation}\label{brKapxy}
[f,g] = \mathop{\sum}\limits_{1\leq i\leq m}
x_iy_i(\del_{x_i}f\cdot\del_{y_i}g + \del_{y_i}f\cdot\del_{x_i}g)
\end{equation}
and
\begin{equation}\label{Qixy}
\begin{alignedat}{2}
&(1+\mathop{\sum}\limits_{1\leq i\leq m}
x_iy_i\del_{x_i}\del_{y_i})f =0&\quad
&\text{for Kap}_{4,0}(2m),\\
&(1+x_1\del_{x_1} + y_1\del_{y_1}+
\mathop{\sum}\limits_{1\leq i\leq m}
x_iy_i\del_{x_i}\del_{y_i})f =0&\quad
&\text{for Kap}_{4,1}(2m).
\end{alignedat}
\end{equation}

For example,
\begin{equation}\label{kapex}
\begin{aligned}
\text{Kap}_{4,0}(2)=\Span(x_1y_1),&\text{Kap}_{4,1}(2)\simeq
\fo_\Pi'(3)\simeq\fvect'(1;(2)),\\
\text{Kap}_{4,0}(4)\simeq\fo_\Pi'(3)\oplus
\fo_\Pi'(3),&\text{Kap}_{4,1}(4)\simeq\text{Kap}_{3}(5)=\fo_\Pi'(5).
\end{aligned}
\end{equation}

\sssec{Gradings and derivations}\label{grader} The commutative subalgebra
$\fh$ in the algebra $\fder(\fg)$ of derivations of the type-2 or type-4
Kaplansky algebra $\fg$, i.e., the subalgebra $\fh$ that determines the
$(\Zee/2)^{2m}$-grading Kaplansky used to construct $\fg$, is not the
maximal torus $\ft$ in $\fder(\fg)$. Clearly, the type-2 or type-4
Kaplansky algebras are $(\Zee/2)^{2m}$-graded by degrees modulo 2 with
respect to each indeterminate $x_i$ and $y_i$; hence,
$\fh=\Span(x_i\del_{x_i}, y_i\del_{y_i}\mid i=1,\dots,m)$. On the other
hand, there exists a $D\in\ft$ commuting with all elements of $\fh$ but
not belonging to $\fh$. Equivalently, there exists a basis of $\fg$
simultaneously homogeneous with respect to the $(\Zee/2)^{2m}$-grading
Kaplansky used and with respect to an extra $\Zee/2$-grading given by
$D$ (which is a second-order operator; see \eqref{Qixy}), and this extra
grading cannot be linearly expressed via the $(\Zee/2)^{2m}$-grading. We
explain why this situation is remarkable.

It might be thought that we should have taken the maximal torus from the
very beginning. The catch is that in all cases we know, except these
Kaplansky algebras, the extra grading operator ``splits" some of the
weight spaces of the previous grading. For each of these Kaplansky
algebras, this is not the case: the weight spaces of the
$(\Zee/2)^{2m}$-grading are already 1-dimensional (except the weight-0
space if we consider the 2-closure of the algebra, but this weight-0
space does not split, anyway). Therefore, the weight spaces cannot be
split further. Hence, it seems there is nowhere the extra grading can
appear from, but it does appear.

We note that the derivation might be given by a differential
operator of order $>1$ but the corresponding grading might still be
``linear" in a sense. We consider the Witt Lie algebra $W_n$ over
$\Kee=\Fee_{2^n}$, where $n>1$. For its basis, we take
$\{e_\alpha\}_{\alpha\in\Kee}$ with the relations $[e_\alpha,
e_\beta]=(\beta-\alpha)e_{\alpha+\beta}$. In fact, $W_n$ is
$\fvect(1; (n))$ over $\Kee$. On $W_n$, there is a natural grading:
$\deg(e_\alpha)=\alpha$.

Now, consider a new grading: $\deg_{new}(e_\alpha)=\alpha^2$, which
resembles the ``nonlinear" gradings of Kaplansky algebras. Indeed,
all weight spaces are 1-dimensional with respect to the old grading,
and the new grading is expressed nonlinearly in terms of the old
grading if $n>1$.

But if the new grading is regarded as $(\Zee/2)^n$-grading (recall
that $\Kee=\Fee_{2^n}=(\Zee/2)^n$ as a vector space), then the new
weight is obtained from the old weight by a linear transformation.
The function $f:\alpha \mapsto \alpha^2$ is linear in the sense that
$f(\alpha+\beta)=f(\alpha)+f(\beta)$, and it is nonlinear in the
sense that it is not true that
\begin{equation}\label{linfun}
f(c\alpha)\neq cf(\alpha)\quad\text{for any }c\in\Kee.
\end{equation}
The condition \eqref{linfun} holds only if $c=0$ or $1$, i.e., for $n=1$.

\paragraph{Gradings not given by derivations}\label{graNOTder}
Because $\Hom(\Zee/q, \Zee/p)=\{0\}$ for primes $q\neq p$, there is
no derivation of the Skryabin algebra $\fby$ that determines its
$\Zee/4$-grading (here $p=3$), see \cite{GL}.

\sssec{The invariant symmetric bilinear forms} Kaplansky also claimed
that each Kaplansky algebra of type 2, 3, or 4 has a nondegenerate
invariant bilinear symmetric form (we call it $K$ here) and several
other interesting properties whose verification ``is quite routine."
Unlike Kaplansky, we think that a lucid proof of these properties is
also of interest. Here, we prove the existence of the invariant form
$K$. The description of $K$ in presence of the alternate form $B$ is
very simple:
\begin{equation}\label{biform}
K(e_u,e_v) = \delta_{u,v}.
\end{equation}
The form $K$ is invariant, i.e.,
\[
K([e_u,e_z],e_v) = K(e_u,[e_z,e_v])
\]
because
\[
\begin{alignedat}{2}
&\text{if }u+z\neq v,&\quad
&\text{then }u\neq z+v,\text{ and both sides vanish, and}\\[1mm]
&\text{if }u+z=v\text{ (and }u=z+v),&\quad
&\text{then the l.h.s.~is }K(B(u,z)e_v,e_v)=B(u,z)\text{ and}\\
&&&\text{the r.h.s.~is }B(z,v)=B(z,u+z)=B(z,u)\\
&&&\text{because }B\text{ is alternate and hence }B(z,z)=0.
\end{alignedat}
\]

We cannot guess how Kaplansky reasoned in the case of the
nonalternate form $B$. In the case of the alternate form $B$, our
argument relies on the invariant form on the Poisson Lie algebra
induced by (the ``desuperization" of) the Berezin
integral\footnote{See \cite{LSh1} for a short summary of the basics
of linear algebra and geometry in a super setting; for a textbook,
see \cite{Lsos} or Bernstein's lectures in \cite{Del}.}
\begin{equation}\label{int}
K(f,g)=\int fg:=\text{the coefficient of the highest term of }fg
\end{equation}
if the Poisson algebra $\fpo_\Pi(n;\uN_s)$ is regarded as a
``desuperization" of the Lie superalgebra $\fpo(0|n)$, i.e., if the
space of $\fpo(0|n)$, the Grassmann superalgebra, is identified with
the algebra of truncated polynomials in even indeterminates.

\ssec{The restricted closures of Kaplansky algebras}\label{KapNotO}
Over $\Fee_2$, the 2-closures of $\fg=\text{Kap}_2(2m)$ and
$\text{Kap}_{4,a}(2m)$, except\footnote{This is a degenerate case:
the algebra is 1-dimensional and its 2-closure is itself.}
$\text{Kap}_{4,0}(2)$, can be described as follows. We set
\begin{equation}\label{2clo}
[\alpha,\beta]=0,\qquad[\alpha,e_u]=\alpha(u)e_u\quad
\text{for any }\alpha,\beta\in V^*,\;e_u\in\fg.
\end{equation}
For a fixed $u\in V$, let $B_u\in V^*$ be the map
\begin{equation}\label{B_u}
B_u:v\mapsto B(u,v)\quad\text{for any }v\in V.
\end{equation}
We can then define squaring by setting
\begin{equation}\label{sq1}
\alpha^{[2]}=\alpha,\qquad e_u^{[2]}=B_u\in V^*.
\end{equation}
The squaring thus defined does indeed satisfy the required
conditions:
\[
[e_u,[e_u, e_v]]=[e_u,B(u,v)e_{u+v}]=
B(u,u+v)B(u,v)e_v=B(u,v)e_v=[B_u, e_v]
\]
and
\[
[\alpha,[\alpha,e_u]]=(\alpha(u))^2e_u=\alpha(u)e_u.
\]

Over an \textit{arbitrary} field $\Kee$ of characteristic 2, the
space of the 2-closure is also $\fg\oplus V^*$, but $\fg$ and $V^*$
are considered over $\Kee$, and squaring is given by the formula
\begin{equation}\label{sq2}
(a\alpha)^{[2]}=a^2\alpha,\qquad
(ae_u)^{[2]}=a^2B_u\in V^*\quad\text{for any }a\in\Kee.
\end{equation}

This description of the 2-closure shows that none of the Lie
algebras $\text{Kap}_{4,a}(2m)$ for $m>2$ is isomorphic to the
simple derived of the orthogonal Lie algebra of the same dimension.
Indeed, the 2-closures of these algebras have different dimensions:
the codimension of the simple derived of the orthogonal algebra in
its 2-closure is much greater than $\dim V^*$. Because $\fo'_I(n)$
is the algebra of zero-diagonal symmetric matrices,
$\dim\fo'_I(n)=\frac12n(n-1)$.

Equation \eqref{dimKap4A} implies that
$\dim\text{Kap}_{4,a}(2m)=\dim\fo'_I(n)$ if $n=2^m+1$ for $a=1$ or
if $n=2^m$ for $a=0$. We therefore wonder if $\text{Kap}_{4,a}(2m)$
is a part of the $\fo'_I(n)$ family. If $n>2$, then the 2-closure of
$\fo'_I(n)$ is the algebra of symmetric traceless matrices, and the
codimension of $\fo'_I(n)$ in its 2-closure is $n-1$ (the dimension
of the space of diagonal matrices of trace 0). And from the above
description, the codimension of $\fo'_I(n)$ in its 2-closure is
$\dim V^*=2m$. Because $n-1>2m$ (if $m>2$), we see that the algebras
$\text{Kap}_{4,a}(2m)$ and $\fo'_I(n)$ are nonisomorphic with the
exceptions $\text{Kap}_{4,0}(2)\simeq \fo'_I(2)$,
$\text{Kap}_{4,1}(2)\simeq \fo'_I(3)$, and
$\text{Kap}_{4,1}(4)\simeq \fo'_I(5)$.

\ssec{General remark on superizations of Lie
algebras}\label{SSgenSuper} The basics of Lie superalgebras for
$p=2$ can be found in \cite{LeP, BGL1}. If $p=2$, there are two
methods which each assigns a simple Lie superalgebra to
\textbf{every} simple Lie algebra (see \cite{BLLS}, where it is
proved that every simple Lie superalgebra is obtained by one of
these two methods from a simple Lie algebra). Here, we apply one of
these methods to Kaplansky algebras. Let $\gr$ be a $\Zee/2$-grading
of $\fg=\fg_\ev\oplus\fg_\od$ and $(\fg,\gr)$ be the minimal
subalgebra of $\overline{\fg}$ containing $\fg$ and all the elements
$x^{[2]}$, where $x\in\fg_\od$. Clearly, there is just one way to
extend the grading $\gr$ to $(\fg,\gr)$. We define squaring by
$x^2:=x^{[2]}$ for any $x\in\fg_\od$ and let $S(\fg,\gr)$ denote the
obtained Lie superalgebra. It is simple.

\sssec{Nonlinear $\Zee/2$-gradings of Kaplansky algebras} The only
known way (until this paper) to obtain a $\Zee/2$-grading on a Lie
algebra amounts to the following. We take an arbitrary linear
function of the weights, more precisely, a homomorphism from the
grading group to $\Zee/2$. Examples of Lie superalgebras
$S(\fg,\gr)$ obtained from these gradings: $\fgl(n)$ produces
$\fgl(k|n-k)$; $\fe(6)$, $\fe(7)$, and $\fe(8)$ produce their
superizations; $\fo_\Pi(2(n+m))$ produces $\fo_{\Pi\Pi}(2n|2m)$ and
$\fpe(2n)$ for $n=m$, whereas $\fh_\Pi(2n;\uN)$ produces
$\fh_\Pi(2k;\tilde\uN|2n-2k)$ and $\fle(n;\tilde\uN)$, if the
coordinates of $\uN=(\tilde\uN,1,\dots,1)$ corresponding to odd
indeterminates are equal to 1 (see \cite{LeP, BGL1}).

The space $V^*$ (more precisely, $\Kee\otimes_{\Fee_2} V^*$, where
$V^*$ is considered over $\Fee_2$) is a torus in the 2-closure of
$\text{Kap}_2(2m)$ or $\text{Kap}_{4,a}(2m)$, and $u\in V$ is
precisely a weight with respect to this torus. That is how we obtain
what we call \textit{linear} superizations of the 2-closures of
$\text{Kap}_{2}(2m)$ and $\text{Kap}_{4,a}(2m)$ (see below).

The Lie algebras $\text{Kap}_{2}(2m)$ give the first (and probably
unique) examples of how to introduce a $\Zee/2$-grading
\textit{nonlinearly}, and there are even two nonequivalent ways to
do this.

Under any superization (linear or not), the even part of the superized
Lie algebra is a Lie subalgebra of the initial Lie algebra. Hence,
there is nothing extraordinary in the fact that the even part of the
superized $\text{Kap}_{2}(2m)\oplus V^*$ is
$\text{Kap}_{4,a}(2m)\oplus V^*$.

The whole $\text{Kap}_{2}(2m)$ cannot enter the even part of the
superized Lie algebra, because the odd part would otherwise be zero.
If $\fg= \text{Kap}_2\oplus V^*$, then $V^*$ cannot be a part of
$\fg_\od$ because $\alpha^{2}=\alpha$ for any $\alpha\in V^*$. Therefore,
$V^*$ must be a part of $\fg_\ev$. Hence, if the whole $\text{Kap}_2$
goes into $\fg_\ev$, there is nothing left for $\fg_\od$.

\sssec{Linear superizations of $\text{Kap}_{2}(2m)$ and
$\text{Kap}_{4,a}(2m)$} Here, we say ``linear" in the sense that
every $e_u$ is homogenous and its parity is a linear function of
$u\in V$ considered over $\Fee_2$.

We define the parity by any element $\varphi\in V^*$ by setting
$p(e_u)=\varphi(u)$. Because the form $B$ is nondegenerate, there
is a unique $v\in V$ such that
\begin{equation}\label{parity}
\varphi=B_v\quad\text{(see \eqref{B_u}), i.e., }\varphi(u)=B(v,u)
\text{ for all }u\in V.
\end{equation}
We let $\varphi_v$ denote this $\varphi$.

To show that two such superizations induced by distinct nonzero
vectors $v$ and $v'$ are isomorphic, it suffices to find a linear
map $M:V\tto V$ such that
\begin{equation}\label{KapSuper}
\begin{alignedat}{2}
&1_2.&\;&M \text{ preserves $B$ for Kap}_{2}(2m),\\
&1_4.&\;&M\text{ preserves $Q$ and hence also preserves $B$
for Kap}_{4,a}(2m),\text{ and}\\
&2.&\;&Mv=v'.\end{alignedat}
\end{equation}
The induced maps
\begin{equation}\label{*}
\tilde M:e_u\mapsto e_{Mu},\qquad
M^*:\varphi\mapsto\varphi\circ M^{-1}\quad
\text{for any }\varphi\in V^*
\end{equation}
then determine an isomorphism between superizations. Indeed, for the
first one,
\[
[\tilde M e_u,\tilde Me_v]=[e_{Mu},e_{Mv}]=
B(Mu,Mv)e_{Mu+Mv}=B(u,v)\tilde M e_{u+v},
\]
and if we also define $P'(e_u)=B(v',u)$, then
\[
P'(e_{Mu})=B(v',Mu)=B(Mv,Mu)=B(v,u)=P(e_u).
\]

\parbegin{Lemma}\label{L1} For $\text{Kap}_{2}(2m)$, an operator $M$
with properties \eqref{KapSuper} exists for any two nonzero vectors
$v$ and $v'$ (we recall that we consider these vectors over $\Fee_2$).
\end{Lemma}

\begin{proof} If $B$ is an alternate bilinear form on a vector space $V$ of dimension $2m$
and $B$ is nondegenerate, then there is an ``alternate basis" for $B$,
i.e., a basis $e^1,\dots,e^{2m}$ of $V$ such that (this is true over
any field of any characteristic; see \cite{Al})
\begin{equation}\label{AltBasis}
B(e^i,e^j)=\begin{cases}\phantom{-}1&\text{if }j=i+m,\\
-1&\text{if }i=j+m,\\
\phantom{-}0&\text{in all other cases,}\end{cases}
\end{equation}
i.e., the Gram matrix of $B$ in this basis is $\begin{pmatrix}
\phantom{-}0_m&1_m\\-1_m&0_m\end{pmatrix}$. \end{proof}

\parbegin{Lemma}\label{L2} Let $B$ and $V$ be as in Lemma \ref{L1}
and $v\in V$ be a nonzero vector. Then there is a basis $e^1,\dots,e^{2m}$
of $V$ satisfying \eqref{AltBasis} such that $e^1 = v$.\end{Lemma}

\begin{proof} We choose any vector $w\in V$ such that $B(v,w)=1$ and
set $e^{m+1}=w$. We set
\[
V_\perp= \{x\in V\mid B(x,v) = B(x,w) = 0\}.
\]
Then $\dim V_\perp = 2m-2$, and the restriction $B_\perp$ of $B$ on
$V_\perp$ is nondegenerate. We choose $e^2,\dots,e^m,e^{m+2},\dots,e^{2m}$
as an alternate basis of $B_\perp$.\end{proof}

Now let $e^1,e^{2},\dots,e^{2m}$ and $\tilde e^1,\tilde e^{2},\dots,\tilde e^{2m}$
be two alternate bases of $V$ such that $e^1=v$ and $\tilde e^1=v'$.
We set $Me^i=\tilde e^i$. Because $B(Me^i,Me^j)=B(\tilde e^i,\tilde e^j)=
B(e^i,e^j)$, it follows that $M$ preserves $B$, and $Mv=v'$.

For $\text{Kap}_{4,a}(2m)$, such an $M$ exists for two nonzero
vectors $v$ and $v'$ (considered over $\Fee_2$) if and only if
$Q(v)=Q(v')$. Hence, there are two linear superizations for each
$\text{Kap}_{4,a}(2m)\oplus V^*$ with the exception of
$\text{Kap}_{4,a}(2)$, where $Q(u)=1$ for any nonzero $u$, which has
only one superization (it is $\fo\fo'_{II}(1|2)$).\footnote{In fact,
the argument with the map \eqref{*} does not prove that the two
superizations of $\text{Kap}_{4,a}(2m)$ are nonisomorphic but only
that there is no isomorphism of the form \eqref{*} between them.
\textbf{Conjecturally}, they are nonisomorphic.}

\sssec{Nonlinear superizations $\text{KapS}_{2,a}(2m)$, and
$\text{KapS}_{4,a,\eps}(2m)$}\label{parKapS} We note that the
superization \eqref{KapS} is \textit{nonlinear}, which means that
the parity is not a linear function of $u$ because it is equal to
$Q(u)+\od$.

Let all spaces defined over $\Fee_2$ be considered over $\Kee$ by
extension of the ground field. We set\footnote{We are not sure which
notation to use here. The $Q_a$ are just \textit{examples} of
quadratic forms with the Arf invariant $a$, while the $Q$ in
\eqref{KapS} can be any quadratic form with the Arf invariant $a$.}
\begin{equation}\label{KapS}
\begin{aligned}
&(\text{KapS}_{2,a}(2m))_\ev:=\text{Kap}_{4,a}(2m)\oplus V^*,\\
&(\text{KapS}_{2,a}(2m))_\od:=\Span(e_u\mid u\in
V,\;u\neq0,\;Q(u)=0)
\end{aligned}
\end{equation}
and define the bracket of even elements with any element and
squaring of the odd elements by eqs.~\eqref{Kapbracket},
\eqref{2clo}, and
\[
(ae_u)^{2}:=(ae_u)^{[2]}=a^2B_u\in V^*\quad\text{(see \eqref{sq2})}.
\]

Let $\text{KapS}_{4,a,\eps}(2m)$ denote the nonlinear superization
of $\text{Kap}_{4,a}(2m)\oplus V^*$ corresponding to a $v\in V$ such
that $Q(v)=\eps$. To describe these Lie superalgebras, we recall the
definition of the parity $\varphi_v$ (see \eqref{parity}),
$\varphi_{v}(u)=B(v,u)$, and consider the following vectors
$v=v_{a,\eps}\in V$ assuming that the quadratic forms $Q_a$ are as
in eq.~\eqref{Q}:
\begin{equation}\label{vepsA}
\begin{aligned}
&v_{0,0}=v_{1,1}=(1,0,\dots,0),\\
&v_{0,1}=(1,0,\dots,0,1,0,\dots,0)\quad\text{(the second 1 is in the
$(m+1)$th position)},\\
&v_{1,0}=(0,1,0,\dots,0)\text{ for }m>1\quad\text{(if $m=1$,
then $Q_1(v)=1$ for any nonzero}\\
&\hskip53mm\text{$v\in V$, and $v_{0,1}$ hence cannot be chosen)}.
\end{aligned}\end{equation}
We set
\begin{equation}\label{KapSreduced}
\begin{aligned}
&\text{KapS}_{4,a,\eps}(2m)_\ev:=\Span(e_u\mid u\neq 0,\,Q_a(u)=1,\,
B(v_{a,\eps}, u)=0) \oplus V^*,\\
&\text{KapS}_{4,a,\eps}(2m)_\od:=\Span(e_u\mid u\neq 0,\,Q_a(u)=1,\,
B(v_{a,\eps}, u)=1).\end{aligned}\end{equation}
(Here, as usual, $B(u,v)=Q(u+v)-Q(u)-Q(v)$, and in this case $Q=Q_a$.)

\paragraph{There are no nonlinear superizations of
$\text{Kap}_{4,a}(2m)$ induced by nonlinear superizations of
$\text{Kap}_{2}(2m)$} In $\text{KapS}_{2,a}(2m)$ corresponding to a
form $Q$, we take the part corresponding to $\text{Kap}_{4,a}(2m)$ with
another form $Q'$. This is a Lie subsuperalgebra. Can we do this? We
can, but this superization fortunately (the classification would
otherwise certainly be a nightmare) coincides with a linear one. This
subsuperalgebra is singled out by the condition $Q'(u)=1$ while its
even part is singled out by this condition together with the extra
condition $Q(u)=1$, which can be replaced with $Q(u)+Q'(u)=0$;
because both $Q$ and $Q'$ should yield the same bilinear form $B$,
the quadratic form $Q+Q'$ degenerates into a linear function.
Therefore, this superization is equivalent to a linear one.

Therefore, up to an isomorphism, there is one linear superization of
$\text{Kap}_{2}(2m)$, this superization\footnote{It would be
interesting to find out if $\text{KapS}_{2}(2m)$ is a deform of a
superization of $\fh_\Pi$. This is clearly not so for
$\fh_\Pi(2k|2m-2k)$ because their dimensions differ (we recall that
$\text{KapS}_{2}(2m)$ contains $V^*$). But it might be a deform of a
larger algebra. \textbf{Conjecturally}, it is not.} is here denoted
by $\text{KapLS}_{2}(2m)$. The three Lie superalgebras
$\text{KapLS}_{2}(2m)$ and $\text{KapS}_{2,a}(2m)$ for $a=0,1$ are
nonisomorphic.

\section{D'inachev\'e}

\ssec{Generalizations of the Jurman construction} We consider
$\fa(2;(g,h))$, the Lie algebra whose space is $\cO(2;(g+h,1))$, and
the bracket of any $F,G\in \cO(2;(g+h,1))$ is given by the formula
(we write $x$ and $y$ to avoid confusion with $p$ and $q$ in the
preceding sections)
\begin{equation}\label{0}
\begin{aligned}
{}[F,G]&=\del_x F \cdot(\del_y + y\del_x^{2^g})G +
(\del_y + y\del_x^{2^g})F \cdot\del_x G\\
&=[F,G]_{P.b.} + y(\del_x F
\cdot \del_x^{2^g}G + \del_x^{2^g} F \cdot\del_x G).
\end{aligned}
\end{equation}
Both $\del_x$ and $\del_y +
y\del_x^{2^g}$ are derivations of $\cO(2;(g+h,1))$ and they
mutually commute and therefore the Jacobi identity holds. (We note that the fact that the
conventional Poisson bracket satisfies the Jacobi identity is a
corollary of the similar properties of $\del_x$ and $\del_y$.)
The first derived $\fa'(2;(g,h))$ of $\fa(2; (g,h))$ is spanned
by all monomials except the highest-degree element $x^{(2^{g+h}-1)}y$.

\sssbegin{Lemma}\label{LemmaIso} We have $\fa'(2;(g,h))/\fc\simeq
\fj(g,h)$ with an isomorphism realized by the expressions
\[
Y_i(0) = x^{(i+1)}y,\qquad Y_i(1) = x^{(i+2)}.
\]
\end{Lemma}

\begin{proof} We directly verify the commutation relations. We first
note that the brackets of $Y_i(0)$ with anything do not contain
additional terms because these terms contain not $\del_y$ but
multiplication by $y$, and $Y_i(0)$ already contains $y$ while
$y\cdot y = 0$. We also note that $[Y_i(1), Y_j(1)]_{P.b.} =0$.
Taking this into account, we see that
%{\footnotesize
\begin{equation}\label{1}
\begin{aligned}
{}[Y_i(0),Y_j(0)]&=x^{(i)}y\cdot x^{(j+1)}+x^{(i+1)}\cdot x^{(j)}y\\
&=\left(\binom{i+j+1}{j+1}+\binom{i+j+1}{i+1}\right)x^{(i+j+1)}y\\
&=\left(\binom{i+j+1}{i+1}+\binom{i+j+1}{j+1}\right)Y_{i+j}(0),
\end{aligned}
\end{equation}
\begin{equation}\label{1.}
\begin{aligned}
{}[Y_i(0),Y_j(1)]&= x^{(i+1)}\cdot xp^{(j+1)}=
\binom{i+j+2}{i+1}x^{(i+j+2)}\\
&=\left(\binom{i+j+1}{i+1}+\binom{i+j+1}{j+1}\right)Y_{i+j}(1).
\end{aligned}
\end{equation}
%}

The statement of eq.~\eqref{1} is clear; that of eq.~\eqref{1.} holds because
if $i+j+1\geq 0$, then
\[
\binom{i+j+1}{i+1}+\binom{i+j+1}{j+1} =
\binom{i+j+1}{i+1}+\binom{i+j+1}{i} = \binom{i+j+2}{i+1},
\]
while if $i+j+1<0$, then $i=j=-1$, and $\binom{i+j+2}{i+1}x^{(i+j+2)}=1$,
i.e., is a constant, which generates the center $\fc$. Therefore, it is
equal to $0$ in the quotient $\fa'(2;(g,h))/\fc$. Hence, in this case,
we also have
$\left(\binom{i+j+1}{i+1}+\binom{i+j+1}{j+1}\right)Y_{i+j}(1)=0$.

We now have
\begin{equation}\label{2}
\begin{aligned}
{}[Y_i(1), Y_j(1)]&= y\left(x^{(i+1)}\cdot x^{(j+1-\eta)} +
x^{(i+1-\eta)}\cdot x^{(j+1)}\right)\\
&= \left(\binom{i+j+2-\eta}{i+1} +
\binom{i+j+2-\eta}{j+1}\right)x^{i+j+2-\eta}y\\
&= \left(\binom{i+j+2-\eta}{i+1} +
\binom{i+j+2-\eta}{j+1}\right)Y_{i+j+1-\eta}(0).
\end{aligned}
\end{equation}
We hence see that the commutation relations are the same
as in $\fj(g,h)$ in all cases.\end{proof}

\ssec{Comparison with known Lie algebras} The direct analog of
bracket \eqref{0} exists in any characteristic $p$ and has the form
\begin{equation}\label{p}
\begin{aligned}
{}[F,G]&= \del_x F \cdot(\del_y +y^{p-1}\del_x^{p^g})G +
(\del_y + y^{p-1}\del_x^{p^g})F \cdot\del_x G \\
&=[F,G]_{P.b.} + y^{p-1}(\del_x F \cdot \del_x^{p^g}G +
\del_x^{p^g} F \cdot\del_x G).
\end{aligned}
\end{equation}
For $p>3$, all finite-dimensional simple Lie algebras are
classified, and this bracket is therefore the bracket of a
known Lie algebra.

\sssbegin{Question}\label{Prob} To which of the filtered deforms
of Lie algebras of Hamiltonian vector fields (see \cite{LeP}) is
the Lie algebra with the bracket \eqref{p} isomorphic?
\end{Question}

\ssec{On further generalizations} We can replace $\del_y +
y\del_x^{2^g}$ with $\del_y + R(y)\del_x^{2^g}$, where $R$ is any
polynomial of a divided degree $\leq \uN(y)$. Conjecturally, the
only $R$ of interest is the monomial of highest possible degree; the
other shapes of $R$ can be reduced to this or a constant. But it
seems that for any $\uN(y)>1$, the result is $\fj(g+N-1,h)$: the
cocycles that make Jurman algebras from $\fh_\Pi'(2; (2,2))$ and
$\fh_\Pi'(2;(3,2))$ change the bracket in precisely this way.

We can consider any number $k$ of
pairs of indeterminates with the bracket
\begin{equation}
\label{3} [F,G] = \mathop{\sum}\limits_{1\leq i\leq k} \del_{x_i} F
\cdot(\del_{y_i} + y_i\del_{x_i}^{2^{g_i}})G + (\del_{y_i} +
y_i\del_{x_i}^{2^{g_i}})F\cdot \del_{x_i} G.
\end{equation}
We note that the $g_i$ can differ for different $i$.

\sssbegin{Lemma}\label{Lbn} The Lie algebra $\fa_\Pi'(2k; (g_1,h_1), \dots,
(g_k,h_k))$ has no center and no homogenous ideals for $k=2$ and
$(g_1,h_1)=(g_2,h_2)=(2,1)$. (Conjecturally, it is
simple.)\end{Lemma}

%$\begin{proof} The lemma is proved using computer-aided study. \end{proof}

\ssec{$\fa_I(2; (g,h))$} The Lie algebra $\fa_I(2; (g,h))$ based on
$\fh_I(2; (g+h, 1))$ can also be generalized in the above way by
beginning with the bracket
\begin{equation}
\label{4} [F,G] = \del_x F \cdot\del_x G + (\del_y + y\del_x^{2^g})F
\cdot(\del_y + y\del_x^{2^g})G
\end{equation}
and generalizing further as indicated above.

\sssbegin{Lemma}\label{Lbn2} The Lie algebra $\fa_I(2k; (g_1,h_1), \dots,
(g_k,h_k))$ has no center and no homogenous ideals for $k=2$ and
$(g_1,h_1)=(g_2,h_2)=(2,1)$. (Conjecturally, it is simple.)\end{Lemma}

%\begin{proof} The lemma is proved using computer-aided study. \end{proof}

{\tiny
\begin{equation}\label{hIdeformed}
\begin{aligned}
&c_{-4}^1 = p\otimes \left (d(p\, q)\wedge d(p{}^{(2)}\, q{}^{(3)})
+d(p\, q{}^{(2)})\wedge d(p{}^{(2)}\, q{}^{(2)}) +d(p\,
q{}^{(3)})\wedge d(p{}^{(2)}\, q)\right )+\dots, \\ %[2mm]
%&+q\otimes \left (d(p\, q)\wedge d(p{}^{(3)}\, q{}^{(2)}) + d(p\,
%q{}^{(2)})\wedge d(p{}^{(3)}\, q)+d(p{}^{(2)}\, q)\wedge
%d(p{}^{(2)}\,q{}^{(2)})\right ) \\
%&+ p{}^{(2)}\otimes \left (d(p\, q)\wedge d(p{}^{(3)}\, q{}^{(3)})
%+d(p)\, \wedge d(p{}^{(3)}\, q{}^{(2)})+(d(p\, q{}^{(3)})\wedge
%d(p{}^{(3)}\, q)\right ) \\
%&+p{}^{(3)}\otimes\left (d(p{}^{(2)}\, q)\wedge d(p{}^{(3)}\,
%q{}^{(3)}) +d(p{}^{(2)}\, q{}^{(2)})\wedge d(p{}^{(3)}\,
%q{}^{(2)})+d(p{}^{(2)}\,q{}^{(3)})\wedge d(p{}^{(3)}\, q)\right )\\
%&+q{}^{(2)}\otimes \left (d(p\, q)\wedge d(p{}^{(3)}\,
%q{}^{(3)})))+d(p\, q{}^{(3)})\wedge d(p{}^{(3)}\, q) +d(p{}^{(2)}\,
%q)\wedge d(p{}^{(2)}\, q{}^{(3)})\right )\\
%& +q{}^{(3)}\otimes \left (d(p\, q{}^{(2)})\wedge d(p{}^{(3)}\,
%q{}^{(3)})+d(p\, q{}^{(3)})\wedge d(p{}^{(3)}\,
%q{}^{(2)})+d(p{}^{(2)}\,
%q{}^{(2)})\wedge d(p{}^{(2)}\, q{}^{(3)})\right )\\
&c_{-4}^2 = p\otimes d(p{}^{(2)}\, q)\wedge d(p{}^{(3)}\, q) +
q\otimes d(p{}^{(3)})\wedge d(p{}^{(3)}\, q) +q\otimes d(p{}^{(2)}\,
q)\wedge d(p{}^{(2)}\, q{}^{(2)}) +\dots,\\%[2mm]
%&+q{}^{(2)}\otimes(d(p{}^{(3)})\wedge d(p{}^{(3)}\, q{}^{(2)}))
%+q{}^{(2)}\otimes(d(p{}^{(2)}\, q)\wedge d(p{}^{(2)}\,
%q{}^{(3)}))+q{}^{(3)}\otimes(d(p{}^{(3)})\wedge
%d(p{}^{(3)}\, q{}^{(3)}))\\
%&+q{}^{(3)}\otimes(d(p{}^{(2)}\, q{}^{(2)})\wedge d(p{}^{(2)}\,
%q{}^{(3)})) +q{}^{(3)}\otimes(d(p{}^{(3)}\, q)\wedge d(p{}^{(3)}\,
%q{}^{(2)})) +(p\, q)\otimes(d(p{}^{(2)}\, q)\wedge d(p{}^{(3)}\,
%q{}^{(2)}))\\
%&+(p\, q)\otimes(d(p{}^{(2)}\, q{}^{(2)})\wedge d(p{}^{(3)}\,
%q))+(p\, q{}^{(2)})\otimes(d(p{}^{(2)}\, q)\wedge d(p{}^{(3)}\,
%q{}^{(3)})) +(p\, q{}^{(2)})\otimes(d(p{}^{(2)}\, q{}^{(3)})\wedge
%d(p{}^{(3)}\, q))\\
%&+p\, q{}^{(3)})\otimes(d(p{}^{(2)}\, q{}^{(2)})\wedge d(p{}^{(3)}\,
%q{}^{(3)}))+(p\, q{}^{(3)})\otimes(d(p{}^{(2)}\,
%q{}^{(3)})\wedge d(p{}^{(3)}\, q{}^{(2)}))\\
&c_{-4}^3=p\otimes d(p\, q)\wedge d(p{}^{(2)}\, q{}^{(3)})
+p\otimes d(p\, q{}^{(2)})\wedge d(p{}^{(2)}\, q{}^{(2)})
+p\otimes d(p\, q{}^{(3)})\wedge d(p{}^{(2)}\, q) +\dots, \\[2mm]
%&q\otimes(d(p\, q)\wedge d(p{}^{(3)}\, q{}^{(2)})) +
%q\otimes(d(p\, q{}^{(2)})\wedge d(p{}^{(3)}\, q)) +
%q\otimes(d(p{}^{(2)}\, q)\wedge d(p{}^{(2)}\, q{}^{(2)})) \\
%&+p{}^{(2)}\otimes(d(p\, q)\wedge d(p{}^{(3)}\, q{}^{(3)}))
%+p{}^{(2)}\otimes(d(p\, q{}^{(2)})\wedge d(p{}^{(3)}\, q{}^{(2)}))
%+p{}^{(2)}\otimes(d(p\, q{}^{(3)})\wedge d(p{}^{(3)}\, q))\\
%& +p{}^{(3)}\otimes(d(p{}^{(2)}\, q)\wedge d(p{}^{(3)}\,
%q{}^{(3)})))+p{}^{(3)}\otimes(d(p{}^{(2)}\, q{}^{(2)})\wedge
%d(p{}^{(3)}\, q{}^{(2)})))+p{}^{(3)}\otimes(d(p{}^{(2)}\,
%q{}^{(3)})\wedge d(p{}^{(3)}\, q))\\
%&+q{}^{(2)}\otimes(d(p\, q)\wedge d(p{}^{(3)}\,
%q{}^{(3)}))+q{}^{(2)}\otimes(d(p\, q{}^{(3)})\wedge d(p{}^{(3)}\,
%q))+q{}^{(2)}\otimes(d(p{}^{(2)}\, q)\wedge
%d(p{}^{(2)}\, q{}^{(3)}))\\
%& +q{}^{(3)}\otimes(d(p\, q{}^{(2)})\wedge d(p{}^{(3)}\,
%q{}^{(3)})))+q{}^{(3)}\otimes(d(p\, q{}^{(3)})\wedge d(p{}^{(3)}\,
%q{}^{(2)})))+q{}^{(3)}\otimes(d(p{}^{(2)}\,
%q{}^{(2)})\wedge d(p{}^{(2)}\, q{}^{(3)}))\\
&c_{-2}^1 = p\otimes d(p{}^{(2)})\wedge d(p{}^{(3)}) +q\otimes
d(p{}^{(2)})\wedge d(p{}^{(2)}\, q) +q{}^{(2)}\otimes
d(p{}^{(2)})\wedge d(p{}^{(2)}\, q{}^{(2)})+\dots,\\%[2mm]
%& q{}^{(3)}\otimes(d(p{}^{(2)})\wedge d(p{}^{(2)}\,
%q{}^{(3)}))+q{}^{(3)}\otimes(d(p{}^{(2)}\, q)\wedge d(p{}^{(2)}\,
%q{}^{(2)}))+(p\, q)\otimes(d(p{}^{(2)})\wedge d(p{}^{(3)}\, q))
%+(p\, q)\otimes(d(p{}^{(3)})\wedge
%d(p{}^{(2)}\, q))\\
%& +(p\, q{}^{(2)})\otimes(d(p{}^{(2)})\wedge d(p{}^{(3)}\,
%q{}^{(2)}))+(p\, q{}^{(2)})\otimes(d(p{}^{(3)})\wedge d(p{}^{(2)}\,
%q{}^{(2)}))+(p\,
%q{}^{(3)})\otimes(d(p{}^{(2)})\wedge d(p{}^{(3)}\, q{}^{(3)}))\\
%& +p\, q{}^{(3)})\otimes(d(p{}^{(3)})\wedge d(p{}^{(2)}\,
%q{}^{(3)}))+(p\, q{}^{(3)})\otimes(d(p{}^{(2)}\, q)\wedge
%d(p{}^{(3)}\, q{}^{(2)}))+(p\,
%q{}^{(3)})\otimes(d(p{}^{(2)}\, q{}^{(2)})\wedge d(p{}^{(3)}\, q))\\
&c_{-2}^2 = p\otimes d(q{}^{(2)})\wedge d(p\, q{}^{(2)}) +q\otimes
d(q{}^{(2)})\wedge d(q{}^{(3)}) +p{}^{(2)}\otimes d(q{}^{(2)})\wedge
d(p{}^{(2)}\, q{}^{(2)}) +\dots,\\%[2mm]
%& +p{}^{(3)}\otimes(d(q{}^{(2)})\wedge d(p{}^{(3)}\,
%q{}^{(2)}))+p{}^{(3)}\otimes(d(p\, q{}^{(2)})\wedge d(p{}^{(2)}\,
%q{}^{(2)}))+(p\,
%q)\otimes(d(q{}^{(2)})\wedge d(p\, q{}^{(3)}))\\
%&+ (p\, q)\otimes(d(q{}^{(3)})\wedge d(p\, q{}^{(2)})) +(p{}^{(2)}\,
%q)\otimes(d(q{}^{(2)})\wedge d(p{}^{(2)}\, q{}^{(3)}))+(p{}^{(2)}\,
%q)\otimes(d(q{}^{(3)})\wedge
%d(p{}^{(2)}\, q{}^{(2)}))\\
%&+(p{}^{(3)}\, q)\otimes(d(q{}^{(2)})\wedge d(p{}^{(3)}\,
%q{}^{(3)})) +(p{}^{(3)}\, q)\otimes(d(q{}^{(3)})\wedge d(p{}^{(3)}\,
%q{}^{(2)}))+(p{}^{(3)}\, q)\otimes(d(p\,
%q{}^{(2)})\wedge d(p{}^{(2)}\, q{}^{(3)}))\\
%&+(p{}^{(3)}\,
%q)\otimes(d(p\, q{}^{(3)})\wedge d(p{}^{(2)}\, q{}^{(2)}))\\[2mm]
&c_{-2}^3 = p\otimes d(p{}^{(2)})\wedge d(p{}^{(3)}) + q\otimes
d(p)\wedge d(p{}^{(3)}\, q) +q\otimes d(p{}^{(2)})\wedge
d(p{}^{(2)}\, q)+\dots,\\%[2mm]
&c_{-2}^4 = p\otimes d(p{}^{(2)})\wedge d(p\, q{}^{(2)}) + p\otimes
d(p{}^{(3)})\wedge d(q{}^{(2)}) +q\otimes d(p{}^{(2)})\wedge
d(q{}^{(3)}) +\dots,\\[2mm]
%&+q\otimes(d(q{}^{(2)})\wedge d(p{}^{(2)}\, q))
%+p{}^{(2)}\otimes(d(p{}^{(2)})\wedge d(p{}^{(2)}\, q{}^{(2)}))
%+p{}^{(3)}\otimes(d(p{}^{(2)})\wedge d(p{}^{(3)}\,
%q{}^{(2)}))\\
%&+p{}^{(3)}\otimes(d(p{}^{(3)})\wedge
%d(p{}^{(2)}\, q{}^{(2)}))+q{}^{(2)}\otimes(d(q{}^{(2)})\wedge
%d(p{}^{(2)}\, q{}^{(2)}))
%+q{}^{(3)}\otimes(d(q{}^{(2)})\wedge d(p{}^{(2)}\,q{}^{(3)}))\\
%&+q{}^{(3)}\otimes(d(q{}^{(3)})\wedge d(p{}^{(2)}\, q{}^{(2)}))+(p\,
%q)\otimes(d(p{}^{(2)})\wedge d(p\, q{}^{(3)})) +(p\, q)\otimes(d(p{}^{(3)})\wedge
%d(q{}^{(3)})) +(p\, q)\otimes(d(q{}^{(2)})\wedge d(p{}^{(3)}\, q))\\
%&+(p\, q)\otimes(d(p\, q{}^{(2)})\wedge d(p{}^{(2)}\,
%q))+(p\, q{}^{(2)})\otimes(d(q{}^{(2)})\wedge
%d(p{}^{(3)}\, q{}^{(2)}))+(p\,
%q{}^{(2)})\otimes(d(p\, q{}^{(2)})\wedge d(p{}^{(2)}\, q{}^{(2)}))\\
%& +(p\, q{}^{(3)})\otimes(d(q{}^{(2)})\wedge d(p{}^{(3)}\,
%q{}^{(3)}))+(p\, q{}^{(3)})\otimes(d(q{}^{(3)})\wedge
%d(p{}^{(3)}\, q{}^{(2)}))+(p\,
%q{}^{(3)})\otimes(d(p\, q{}^{(2)})\wedge d(p{}^{(2)}\, q{}^{(3)}))\\
%& +(p\, q{}^{(3)})\otimes(d(p\, q{}^{(3)})\wedge d(p{}^{(2)}\,
%q{}^{(2)})))+(p{}^{(2)}\, q)\otimes(d(p{}^{(2)})\wedge
%d(p{}^{(2)}\, q{}^{(3)}))+(p{}^{(2)}\,
%q)\otimes(d(p{}^{(2)}\, q)\wedge d(p{}^{(2)}\, q{}^{(2)}))\\
%& +(p{}^{(3)}\, q)\otimes(d(p{}^{(2)})\wedge d(p{}^{(3)}\,
%q{}^{(3)}))+(p{}^{(3)}\, q)\otimes(d(p{}^{(3)})\wedge
%d(p{}^{(2)}\, q{}^{(3)}))+(p{}^{(3)}\,
%q)\otimes(d(p{}^{(2)}\, q)\wedge d(p{}^{(3)}\, q{}^{(2)}))\\
%& +(p{}^{(3)}\, q)\otimes(d(p{}^{(2)}\, q{}^{(2)})\wedge
%d(p{}^{(3)}\, q))
&c_0 = p\otimes d(q)\wedge d(p\, q) +p{}^{(2)}\otimes d(q)\wedge
d(p{}^{(2)}\, q)) +p{}^{(3)}\otimes d(q)\wedge d(p{}^{(3)}\, q) +\dots,\\[2mm]
%&+p{}^{(3)}\otimes(d(p\, q)\wedge d(p{}^{(2)}\, q))
%+q{}^{(3)}\otimes(d(p\, q)\wedge d(p\, q{}^{(2)})) +(p\,
%q{}^{(3)})\otimes(d(p\, q)\wedge d(p{}^{(2)}\,q{}^{(2)}))\\
%&+(p\, q{}^{(3)})\otimes(d(p\,
%q{}^{(2)})\wedge d(p{}^{(2)}\, q))+(p{}^{(2)}\,
%q{}^{(3)})\otimes(d(p\, q)\wedge d(p{}^{(3)}\, q{}^{(2)}))
%+(p{}^{(2)}\, q{}^{(3)})\otimes(d(p\, q{}^{(2)})\wedge d(p{}^{(3)}\, q))\\
%&+(p{}^{(3)}\,q{}^{(3)})\otimes(d(p{}^{(2)}\, q)\wedge d(p{}^{(3)}\, q{}^{(2)}))
%+(p{}^{(3)}\, q{}^{(3)})\otimes(d(p{}^{(2)}\, q{}^{(2)})\wedge d(p{}^{(3)}\, q))\\[2mm]
&c_{2}^1 = q{}^{(3)}\otimes d(q)\wedge d(p{}^{(2)}) )+p\,
q{}^{(3)} \otimes d(q)\wedge d(p{}^{(3)}) + p\,
q{}^{(3)} \otimes d(p{}^{(2)})\wedge d(p\, q) +\dots,\\%[2mm]
%&+ p{}^{(2)}\, q{}^{(3)} \otimes d(p{}^{(2)})\wedge d(p{}^{(2)}\,
% q) + p{}^{(3)}\, q{}^{(3)} \otimes d(p{}^{(2)})\wedge
% d(p{}^{(3)}\, q) +p{}^{(3)}\,
% q{}^{(3)}\otimes d(p{}^{(3)})\wedge d(p{}^{(2)}\, q) \\[2mm]
&c_{2}^2 = p{}^{(3)}\otimes d(p)\wedge d(q{}^{(2)}) + p{}^{(3)}\, q
\otimes d(p)\wedge d(q{}^{(3)}) + p{}^{(3)}\, q \otimes
d(q{}^{(2)})\wedge d(p\, q)+\dots,\\
%& + p{}^{(3)}\, q{}^{(2)} \otimes d(q{}^{(2)})\wedge d(p\,
%q{}^{(2)}) + p{}^{(3)}\, q{}^{(3)} \otimes d(q{}^{(2)})\wedge d(p\,
%q{}^{(3)}) + p{}^{(3)}\, q{}^{(3)} \otimes d(q{}^{(3)})\wedge d(p\,
%q{}^{(2)}) \\[2mm]
&c_{2}^3 = q{}^{(3)}\otimes d(q)\wedge d(q{}^{(2)}) + p\, q{}^{(3)}
\otimes d(q)\wedge d(p\, q{}^{(2)}) + p\, q{}^{(3)} \otimes
d(q{}^{(2)})\wedge d(p\, q)+\dots,\\
%& + p{}^{(2)}\, q{}^{(3)} \otimes d(q)\wedge d(p{}^{(2)}\,
%q{}^{(2)}) + p{}^{(2)}\, q{}^{(3)} \otimes d(q{}^{(2)})\wedge
%d(p{}^{(2)}\, q) + p{}^{(3)}\, q{}^{(3)} \otimes d(q)\wedge
%d(p{}^{(3)}\,q{}^{(2)}) \\
%&+ p{}^{(3)}\, q{}^{(3)} \otimes d(q{}^{(2)})\wedge d(p{}^{(3)}\,
%q)) + p{}^{(3)}\, q{}^{(3)} \otimes d(p\, q)\wedge d(p{}^{(2)}\,
%q{}^{(2)}) + p{}^{(3)}\, q{}^{(3)} \otimes d(p\,
%q{}^{(2)})\wedge d(p{}^{(2)}\, q) \\[2mm]
&c_{2}^4 = p{}^{(3)}\otimes d(p)\wedge d(p{}^{(2)}) + p{}^{(3)}\, q
\otimes d(p)\wedge d(p{}^{(2)}\, q) + p{}^{(3)}\, q \otimes
d(p{}^{(2)})\wedge d(p\, q) +\dots,\\[2mm]
%& +p{}^{(3)}\, q{}^{(2)} \otimes d(p)\wedge d(p{}^{(2)}\, q{}^{(2)})
%+ p{}^{(3)}\, q{}^{(2)} \otimes d(p{}^{(2)})\wedge d(p\, q{}^{(2)})
%+ p{}^{(3)}\, q{}^{(3)} \otimes d(p)\wedge d(p{}^{(2)}\,
%q{}^{(3)}) \\
%&+p{}^{(3)}\,
% q{}^{(3)} \otimes d(p{}^{(2)})\wedge d(p\, q{}^{(3)}) + p{}^{(3)}\, q{}^{(3)} \otimes(d(p\, q)\wedge d(p{}^{(2)}\,
% q{}^{(2)}) + p{}^{(3)}\, q{}^{(3)} \otimes d(p\,
%q{}^{(2)})\wedge d(p{}^{(2)}\, q) \\[2mm]
&c_{6} = p{}^{(3)}\, q{}^{(3)} \otimes d(p)\wedge d(q).
\end{aligned}
\end{equation}}

\parbegin{Lemma}\label{nonlin} For $\fg:=\fh_I(2;(2,2))$,
each cocycle \eqref{hIdeformed} representing the weight elements of
$H^2(\fg; \fg)$ is integrable, and all except $c_{-2}^3$ (see
eq.~\eqref{c-23}) are \textbf{linearly} integrable.
\end{Lemma}

\begin{proof} Computer-aided. The nonlinear deform (cf.~\cite{BLW}) is
\[
[\cdot ,\cdot]_{\hbar}=[,]+c_{-2}^3 \hbar+A \hbar^2+B\hbar^3,
\]
where {\tiny
\begin{equation}\label{c-23}
\begin{aligned}
A={}&p \, q {}^{(3)}\otimes
d(p {}^{(2)}\, q {}^{(2)})\wedge d(p {}^{(3)}\,
q {}^{(3)}) +p \, q {}^{(3)}\otimes d(p {}^{(2)}\,
q {}^{(3)})\wedge d(p {}^{(3)}\, q {}^{(2)})\\
&{}+p {}^{(2)}\,
q {}^{(2)}\otimes d(p {}^{(3)}\, q )\wedge d(p {}^{(3)}\, q {}^{(3)})+q {}^{(3)}\otimes d(p \, q {}^{(2)})\wedge d(p {}^{(3)}\,
q {}^{(3)})\\
&{}+q {}^{(3)}\otimes d(p \, q {}^{(3)})\wedge
d(p {}^{(3)}\, q {}^{(2)})+p \,
q {}^{(2)}\otimes d(p {}^{(2)}\, q )\wedge d(p {}^{(3)}\, q {}^{(3)})\\
&{}+p {}^{(2)}\, q \otimes d(p {}^{(3)}\, q )\wedge d(p {}^{(3)}\,
q {}^{(2)}) +q {}^{(2)}\otimes d(p \, q )\wedge d(p {}^{(3)}\,
q {}^{(3)})\\
&{}+q \otimes d(p {}^{(3)})\wedge d(p {}^{(3)}\,
q )+q \otimes d(p \, q )\wedge d(p {}^{(3)}\,
q {}^{(2)}),\\[2mm]
B={}&q {}^{(3)}\otimes d(p {}^{(3)}\, q {}^{(2)})\wedge
d(p{}^{(3)}\, q{}^{(3)}). \hskip 3cm\qed
\end{aligned}
\end{equation}}
\noqed\end{proof}

\paragraph{Claim: The Lie algebra $\fa_I(2; (g,h))$ is a
deform of $\fh_I(2; (g+h, 1))$}\label{Claim} To prove this for the
smallest values of $(g,h)$, we list all infinitesimal deforms of
$\fh_I(2;(2,2))$. For the cochain
$F\otimes\left(dG_1\wedge\cdots\wedge dG_n\right)$, where $F,
G_1,\dots,G_n\in \fh_I(2;(g+h,1))$, its weight is equal to
\begin{equation}\label{innerdeg}
((\deg_p(F) - \mathop{\sum}\limits_{1\leq i\leq n}
\deg_p(G_i)) \bmod 2,\quad (\deg_q(F) - \mathop{\sum}\limits_{1\leq
i\leq n} \deg_q(G_i)) \bmod 2).
\end{equation}
We note that this grading is induced by elements of a maximal torus,
more specifically, by $p^{(2)}$ and $q^{(2)}$. This grading is
therefore modulo 2 and is not a $\Zee$-grading. This algebra also has
the outer grading $\deg_{\mathrm{out}}$ given by
\begin{equation} \label{outdeg}
\deg(p)=\deg(q)=1,\quad
\deg_{\mathrm{out}}(F)=\deg(F)-2,\quad
\deg_{\mathrm{out}}(dF)=2-\deg(f).
\end{equation}
The cocycles \eqref{gh=31cocycles} are all of weight $\{0,0\}$. They
are indexed in accordance with $\deg_{\mathrm{out}}$.

\paragraph{Question}\label{gh=31prob} How to interpret the non-Jurman
cocycles \`a la Proposition \ref{Jcocycles} for the other values of
$(g, h)$? For example, for $(g, h)=(3,1)$ and $(2,2)$, i.e., for the
deformations of $\fh_\Pi'(2;(3,2))\simeq\fh_\Pi'(2;(2,3))$, the
Jurman cocycle deforming $\fh_\Pi'(2;(3,2))$ into $\fj(3,1)$ is
$c_{-2,8}$, and the Jurman cocycle deforming $\fh_\Pi'(2;(2,3))$
into $\fj(2,2)$ is $c_{4,-2}$, see \eqref{gh=31cocycles}.
{\tiny
\begin{equation}\label{gh=31cocycles}
\begin{aligned}
&c_{0,-8}=p\otimes d\left(p\, q{}^{(4)}\right)\wedge
d\left(p\, q{}^{(5)}\right) + p\otimes
d\left(q{}^{(5)}\right)\wedge d\left(p{}^{(2)}\,
q{}^{(4)}\right) +q\otimes d\left(p\,
q{}^{(4)}\right)\wedge d\left(q{}^{(6)}\right)+\dots,\\
&c_{1,-7}=p\otimes d\left(q{}^{(4)}\right)\wedge d\left(p\,
q{}^{(4)}\right) +q\otimes d\left(q{}^{(4)}\right)\wedge
d\left(q{}^{(5)}\right)+p{}^{(2)}\otimes
d\left(q{}^{(4)}\right)\wedge d\left(p{}^{(2)}\,
q{}^{(4)}\right)+\dots,\\
&c_{4,-4}=p{}^{(3)}\otimes d\left(q\right)\wedge
d\left(q{}^{(4)}\right) + p{}^{(3)}\, q \otimes
d\left(q\right)\wedge d\left(q{}^{(5)}\right) + p{}^{(3)}\, q
\otimes d\left(q{}^{(2)}\right)\wedge
d\left(q{}^{(4)}\right) +\dots,\\[2mm]
%& p{}^{(3)}\,q{}^{(2)} \otimes d\left(q\right)\wedge
%d\left(q{}^{(6)}\right) + p{}^{(3)}\, q{}^{(2)} \otimes
%d\left(q{}^{(3)}\right)\wedge d\left(q{}^{(4)}\right) + p{}^{(3)}\,
%q{}^{(3)} \otimes d\left(q\right)\wedge
%d\left(q{}^{(7)}\right) +\\
%& p{}^{(3)}\, q{}^{(3)} \otimes d\left(q{}^{(2)}\right)\wedge
%d\left(q{}^{(6)}\right) + p{}^{(3)}\, q{}^{(3)} \otimes
%d\left(q{}^{(3)}\right)\wedge d\left(q{}^{(5)}\right) + p{}^{(3)}\,
%q{}^{(4)} \otimes d\left(q{}^{(4)}\right)\wedge
%d\left(q{}^{(5)}\right) +\\
%& p{}^{(3)}\, q{}^{(5)} \otimes d\left(q{}^{(4)}\right)\wedge
%d\left(q{}^{(6)}\right) + p{}^{(3)}\, q{}^{(6)} \otimes
%d\left(q{}^{(4)}\right)\wedge d\left(q{}^{(7)}\right) + p{}^{(3)}\,
%q{}^{(6)} \otimes
%d\left(q{}^{(5)}\right)\wedge d\left(q{}^{(6)}\right) \\[2mm]
&c_{4,-2}=p{}^{(3)}\otimes d\left(q\right)\wedge
d\left(q{}^{(2)}\right) + p{}^{(3)}\, q \otimes
d\left(q\right)\wedge
d\left(q{}^{(3)}\right) + p{}^{(3)}\, q{}^{(2)} \otimes
d\left(q\right)\wedge
d\left(q{}^{(4)}\right) +\dots,\\
%& p{}^{(3)}\, q{}^{(2)} \otimes d\left(q{}^{(2)}\right)\wedge
%d\left(q{}^{(3)}\right) + p{}^{(3)}\, q{}^{(3)} \otimes
%d\left(q\right)\wedge d\left(q{}^{(5)}\right) + p{}^{(3)}\,
%q{}^{(4)} \otimes
%d\left(q\right)\wedge d\left(q{}^{(6)}\right) +\\
%& p{}^{(3)}\, q{}^{(4)} \otimes d\left(q{}^{(2)}\right)\wedge
%d\left(q{}^{(5)}\right) + p{}^{(3)}\, q{}^{(5)} \otimes
%d\left(q\right)\wedge d\left(q{}^{(7)}\right) + p{}^{(3)}\,
%{}^{(5)} \otimes d\left(q{}^{(3)}\right)\wedge
%d\left(q{}^{(5)}\right) +\\
%& p{}^{(3)}\, q{}^{(6)} \otimes d\left(q{}^{(2)}\right)\wedge
%d\left(q{}^{(7)}\right) + p{}^{(3)}\, q{}^{(6)} \otimes
%d\left(q{}^{(3)}\right)\wedge d\left(q{}^{(6)}\right)+p{}^{(3)}\,
%q{}^{(6)} \otimes d\left(q{}^{(4)}\right)
%\wedge d\left(q{}^{(5)}\right) \\[2mm]
&c_{1,-5}=p\otimes d(q{}^{(2)})\wedge d(p\, q{}^{(4)})
+p\otimes d(p\, q{}^{(2)})\wedge d(q{}^{(4)}) +\dots,\\[2mm]
&c_{0,-4}=p\otimes d(p\, q{}^{(2)})\wedge d(p\, q{}^{(3)}) +p\otimes
d(q{}^{(3)})\wedge d(p{}^{(2)}\, q{}^{(2)}) +\dots,\\
&c_{-1,-5}=p\otimes d(p{}^{(2)})\wedge d(p\, q{}^{(6)})
+p\otimes d(p{}^{(3)})\wedge d(q{}^{(6)})+\dots,\\
&c_{-2,-6}=p\otimes d(p\, q{}^{(4)})\wedge d(p{}^{(3)}\, q{}^{(3)})
+ q\otimes d(p\, q{}^{(4)})\wedge d(p{}^{(2)}\, q{}^{(4)}) +\dots,\\[2mm]
&c_{-2,-4}=p\otimes d(p\, q{}^{(2)})\wedge d(p{}^{(3)}\, q{}^{(3)})
+p\otimes d(p{}^{(3)}\, q)\wedge d(p\, q{}^{(4)}) +\dots,\\
&c_{-1,-3}=p\otimes d(q{}^{(2)})\wedge d(p{}^{(3)}\, q{}^{(2)})
+p\otimes d(p{}^{(2)}\, q)\wedge d(p\, q{}^{(3)}) +\dots,\\
&c_{0,-2}=p\otimes d(p\, q)\wedge d(p\, q{}^{(2)}) +p\otimes
d(q{}^{(2)})\wedge d(p{}^{(2)}\, q) +\dots,\\
&c_{2, 0}=p{}^{(2)}\otimes d(p)\wedge d(q) + p\,
q{}^{(2)} \otimes d(q)\wedge d(q{}^{(2)}) +\dots,\\[2mm]
&c_{-2,-2}=p\otimes d(p\, q{}^{(2)})\wedge d(p{}^{(3)}\, q)
+q\otimes d(q)\wedge d(p{}^{(3)}\, q{}^{(3)}) +\dots,\\[2mm]
&c_{-2,0}=p\otimes d(p{}^{(2)})\wedge d(p{}^{(2)}\, q)
+p\otimes d(p\, q)\wedge d(p{}^{(3)}) +\dots,\\
&c_{-4,-2}=p\otimes d(p{}^{(3)})\wedge d(p{}^{(3)}\, q{}^{(3)}) +
q\otimes d(p{}^{(3)})\wedge d(p{}^{(2)}\, q{}^{(4)}) +\dots,\\[2mm]
&c_{-4,0}=p\otimes d(p{}^{(3)})\wedge d(p{}^{(3)}\, q)
+q\otimes d(p{}^{(3)})\wedge d(p{}^{(2)}\, q{}^{(2)}) +\dots,\\
&c_{0,4}=(q{}^{(4)}\otimes \left(d(p)\wedge d(q)\right) +
(p{}^{(2)}\,
q{}^{(3)})\otimes d(p)\wedge d(p{}^{(2)}) +\dots,\\[2mm]
&c_{0,6}=q{}^{(6)}\otimes d\left(p\right)\wedge d\left(q\right) +
p{}^{(2)}\,
q{}^{(5)} \otimes d\left(p\right)\wedge d\left(p{}^{(2)}\right)+\dots,\\
%& + p\, q{}^{(7)} \otimes d\left(p\right)\wedge d\left(p\,
%q{}^{(2)}\right) + p{}^{(3)}\, q{}^{(6)} \otimes
% d\left(p\right)\wedge d\left(p{}^{(3)}\, q\right) \\
%& + p{}^{(2)}\, q{}^{(7)} \otimes d\left(q\right)\wedge
%d\left(p{}^{(3)}\, q\right) + p{}^{(2)}\, q{}^{(7)} \otimes
%d\left(p{}^{(2)}\right)\wedge
%d\left(p\, q{}^{(2)}\right) \\[2mm]
&c_{-2,8}=q{}^{(7)}\otimes d(p)\wedge d(p{}^{(2)})
+p\,q{}^{(7)}\otimes d(p)\wedge d(p{}^{(3)}) + p{}^{(2)}\,q{}^{(7)}
\otimes d(p{}^{(2)})
\wedge d(p{}^{(3)}).
\end{aligned}
\end{equation}}

\sssbegin{Conjecture}\label{conj1} The Lie algebra
\textit{Kap}${}_{4,B}(2m)$ is not isomorphic to $\fpo_\Pi(2m;\uN_s)$,
and \textit{Kap}${}_{2}(2m)$ is not isomorphic to $\fh_\Pi(2m;\uN)$.
\end{Conjecture}

We verified this for small $m$. For $m=1$, $\text{Kap}_{4,B}(2)$ is
isomorphic to $\fo'(3)\oplus\fc$, where $\fc$ is the 1-dimensional
trivial center and is hence not isomorphic to $\fpo_\Pi(2;\uN_s)$,
which is solvable. For $m=2$, computer-aided computations show that
the infinitesimal deformation corresponding to \eqref{brKap-deformed}
is a nontrivial cocycle. To prove the conjecture, we must show that
the cocycle is also not semitrivial. Of course, what we really need
to know is what $\text{Kap}_{4,B}(2m)$ and its subalgebras
$\text{Kap}_{4,a}(2m)$ \textit{are} isomorphic to. We present some
plausible conjectures.

\sssbegin{Conjectures}\label{conj2}
1. The Lie algebra \textit{Kap}${}_{4,1}(2m)$ is a deform of the
subalgebra in the Poisson algebra $\fpo(2m;\uN_s)$ generated by
functions $f\in\cO(2m;\uN_s)$ satisfying
$\mathop{\sum}\limits_{1\leq i\leq3}\frac{\del^2 f}{\del p_i\del q_i}=0$.
(The quotient of this subalgebra modulo center is isomorphic
to $\fs\fl\fh(2m)$; see \cite{LeP}.)

2. The Lie algebra \textit{Kap}${}_{4,1}(2m)$ is a deform of
$\fo_I'(2m+1;\uN_s)$ while \textit{Kap}${}_{4,0}(2m)$ is a deform of
a subalgebra in $\fo_I'(2m;\uN_s)$ (see \cite{LeP}).
\end{Conjectures}

The dimension of $H^2(\fg;\fg)$ is big and grows quickly with $m$.
How can we select the needed deform? The Poisson algebra and its
subalgebra consisting of harmonic functions have a center generated
by constants, while $\text{Kap}_{4,1}(2m)$ is simple. Therefore, in
the huge space of cocycles representing infinitesimal deformations,
we need only select cocycles of the form
\begin{equation}\label{nonnter}
f\otimes d(1)\wedge d(g)+\dots
\end{equation}
and compare the global deforms corresponding to such cocycles with
$\text{Kap}_{4,1}(2m)$. For small $m$, $\dim H^2(\fg;\fg)$ does not
explode yet. For $m=2$ and $m=3$, we have $\dim H^2(\fg;\fg)=34$; all
cocycles are integrable and all global deforms corresponding to them
(if a representative is chosen carefully by means of coboundaries)
are linear in the deformation parameter. For $m=2$ and $m=3$, there
is only one cocycle of the form \eqref{nonnter} (up to coboundaries).
These cocycles are of degree 2. In degree 2, there is only one
cocycle for $m=3$, and there are five cocycles for $m=2$. Further
investigations show that Conjecture 1 only holds for $m=2$; for
$m=3$, the two algebras to be compared have different numbers of
central extensions.

\ssec{How to establish nonisomorphicy?}\label{sssIso} Skryabin
\cite{Sk} classified the filtered deforms of Hamiltonian Lie
algebras $\fh_\Pi(2m;\uN)$. It remains to select which of them is
the simple Lie algebra
$\text{Kap}_{4,B}(2m)/\fc\simeq\text{Kap}_{2}(2m)$. We have not yet
performed such an identification.

To find out if two given Lie algebras of the same dimension are
isomorphic, Eick considered the following invariants in
\cite{Ei}:\footnote{An almost exact quotation from \cite{Ei}: ``We
say that a derivation $d\in\fder(\fg)$ is \textit{$p$-nilpotent} if
$d^p=0$ holds. For a $p$-nilpotent derivation $d$, we define its
exponential matrix $\exp d:=\mathop{\sum}\limits_{0\leq i\leq
p-1}\frac{d^i}{i!}$. We call a $p$-nilpotent derivation $d$ an
\textit{annihilator} if $d^i(X) d^j(Y)=0$ for all $X,Y\in\fg$ and
$i,j\geq0$ with $i+j\geq p$. Let $\Ann(\fg)\subset \fder (\fg)$ denote
the subset of annihilators. We define $\Exp(\fg)$ to be the subgroup
of $\Aut(\fg)$ generated by $\{\exp(d)\mid d \in\Ann(\fg)\}$. We note
that the order of every element $\exp(d)$ is equal to either $p$ or
$1$. Hence, $\Exp(\fg)$ is a subgroup of $\Aut(\fg)$ generated by
automorphisms of order $p$."} $\dim H^1(\fg; \fg)$ or rather $\dim
\fder (\fg)$, the order of the group $\Aut(\fg)$, the number of
elements in $\Ann(\fg)$, and the order of $\Exp(\fg)$.

Speaking of deforms, we can consider the action of $\Aut(\fg)$ on
the space of infinitesimal deformations, as in \cite{KCh, Ch}.

For algebras of small dimension, there is still another approach,
at least theoretically. We can compare identities that the algebras
satisfy. A.~A.~Kirillov formulated the following analog of the
Amitsur--Levitzki theorem, whose proof was only preprinted in the
Keldysh Institute of Applied Mathematics in the 1980s (see \cite{KOU}
for a translation of one such preprint; the other preprints with related
results by Kirillov, Kontsevich, and Molev have not yet been translated,
but they were at least reviewed by Molev).

\sssbegin{Theorem}[\cite{Ki}] Let $\fg$ be a simple Lie algebra of
vector fields over a field of characteristic $0$. Let
\begin{equation}\label{T_k}
a_k(X_1, \dots, X_k)=\mathop{\sum}\limits_{\sigma\in
S_k}(-1)^{\sign\sigma}\ad_{X_{\sigma(1)}} \dots \ad_{X_{\sigma(k)}}.
\end{equation}
The identity $a_k(X_1,\dots,X_k)\equiv0$ for any
$X_1,\dots,X_k \in\fg$ holds

\textrm{a.} for $k\geq (n+1)^2$ if $\fg=\fvect(n)$,

\textrm{b.} for $k\geq n(2n+5)$ if $\fg=\fh(2n)$, and

\textrm{c.} for $k\geq 2n^2+5n+5$ if $\fg=\fk(2n+1)$.
\end{Theorem}

Dzhumadildaev suggested an interesting modification of emphasis in
this train of thought, finding a hidden supersymmetry for an analog
of antisymmetrizors with just $x$ instead of $\ad_{x}$ in
\eqref{T_k}. He also showed a relation to strongly homotopy algebras
(for further details, see \cite{Dzhu} and \cite{LL}).

\section*{Acknowledgements}
We are thankful to P.~Grozman for his wonderful package
\textit{SuperLie} (see \cite{Gr}). Special thanks are due to the
referee for the helpful comments and to P.~Zusmanovich, who informed
us about important papers \cite{Ei} and \cite{SkT1}. S.B.~was
supported in part by the grant AD~065~NYUAD.

\bibliographystyle{mrl}

\end{document}